\documentclass[11pt]{article}
\usepackage{IEEEtrantools}
\usepackage{mathtools, nccmath}
\usepackage{latexsym}
\usepackage{amsfonts}
\usepackage{amssymb}
\usepackage{amsmath,amsfonts,amsthm}
\usepackage{mathrsfs}
\usepackage{enumerate}
\usepackage{dsfont}
\usepackage{hyperref}
\usepackage{graphicx}
\usepackage{epstopdf}
\usepackage{tikz}

\usetikzlibrary{shapes,arrows.meta,positioning}
\usepackage{cleveref}
\usepackage{tikz-cd}
\usepackage{psfrag}

\usepackage{multirow}
\usepackage{tikz}
\usetikzlibrary{arrows,decorations.markings}
\usepackage{bookmark}
\usepackage{hyperref}
\hypersetup{
     colorlinks   = true,
     citecolor    = blue
}

\newcommand{\newsection}[1]{\setcounter{equation}{0}
\setcounter{dfn}{0}
\section{#1}}

\newtheorem{dfn}{Definition}[section]
\newtheorem{thm}[dfn]{Theorem}
\newtheorem{lmma}[dfn]{Lemma}
\newtheorem{ppsn}[dfn]{Proposition}
\newtheorem{crlre}[dfn]{Corollary}

\newtheorem{rmrk}[dfn]{Remark}
\newtheorem{notation}[dfn]{Notation}

\DeclareMathOperator*{\dprime}{\prime \prime}

\newcommand{\bbc}{\mathbb{C}}
\newcommand{\bbz}{\mathbb{Z}}
\newcommand{\bbn}{\mathbb{N}}

\def \qed { \mbox{}\hfill
$\Box$\vspace{1ex}}

\title{Sections and Chapters}

\begin{document}

\tikzset{->-/.style={decoration={
  markings,
  mark=at position #1 with {\arrow{>}}},postaction={decorate}}}
  \tikzset{-<-/.style={decoration={
  markings,
  mark=at position #1 with {\arrow{<}}},postaction={decorate}}}

\author{\sc{Keshab Chandra Bakshi\,,\,Satyajit Guin\,,\,Guruprasad}}
\title{Conjugate pairs of Hadamard subfactors and vertex models}
\maketitle


\begin{abstract}
We show that any two Hadamard subfactors arising from a pair of distinct complex Hadamard matrices of order $3$ are either equal or conjugate by a unitary in the relative commutant of their intersection. Moreover, when the Hadamard subfactors are not equal, we prove the factoriality of their intersection, and it turns out to be a vertex model subfactor. We compute the first relative commutant and characterize this subfactor by identifying it with a particular type of Krishnan-Sunder subfactor. A few key invariants, including the Pimsner-Popa probabilistic number, the angle, and the Connes-St{\o}rmer relative entropy for the pair of Hadamard subfactors are computed to understand their relative position.
\end{abstract}
\bigskip

{\bf AMS Subject Classification No.:} {\large 46}L{\large 37}, {\large 46}L{\large 55}, {\large 46}L{\large 10}, {\large 37}A{\large 35}.
\smallskip


\hypersetup{linkcolor=blue}
\bigskip



\newsection{Introduction}\label{Sec 1}

Jones' subfactor theory may be thought of as a quantization of a closed subspace in a Hilbert space \cite{Jo, Jo2}. The study of multiple subfactors was initiated by Ocneanu. He proposed the concept of ‘maximal atlas’ for a compatible family of ‘finite-index’ bimodules arising from the subfactors \cite{Oc}. Ocneanu’s work leads Jones to propose the study of `two subfactors' as a quantized version of a pair of subspaces in a Hilbert space \cite{Jo1}. A few important invariants for two subfactors are the Pimsner-Popa probabilistic number \cite{PP}, the Sano-Watatani angle operator \cite{SW}, the interior and exterior angles \cite{BDLR}, and the Connes-St{\o}rmer relative entropy \cite{CS}. In practice, however, computations of these invariants are often challenging. In the absence of a general theory for two subfactors, and to gain fruitful insight, it is imperative that we first investigate some special classes of two subfactors in order to hope for building a general theory. Motivated by this goal, in \cite{BG1,BG2} the first two authors initiated an investigation of a special class of two subfactors, namely a pair of spin model subfactors, also known as the Hadamard subfactors. Importance of this class of subfactors has been emphasized by Jones \cite{Jo2}; however, not much is known about the structure of these subfactors in general.

Recent investigation carried out in \cite{BG1,BG2} helps us gain insights about pairs of Hadamard subfactors, and provide a starting point for further investigations of two subfactor theory. Several interesting results are obtained, some in the general situation and some in the concrete cases, including explicit computation of the invariants. Hadamard subfactors arising from a pair of complex Hadamard matrices of order $2$ (see \cite{BG1}), as well as Hadamard inequivalent complex Hadamard matrices of order $4$ (see \cite{BG2}) have been investigated in depth. It has been observed that there are sharp contrasts between these two cases. One of the major contrasts is that in the $2\times 2$ case, the intersection is a factor of fixed index (equal to $4$), while in the $4\times 4$ situation, the value of the index lies in the set $\{4n:n\geq 2\}\cup\{\infty\}$ including all possibilities. On the other hand, the intersection in the $2\times 2$ situation is non-irreducible, while that in the $4\times 4$ situation is irreducible in some cases. Moreover, the relative entropy in the $2\times 2$ situation depends on the input matrices, while that in the $4\times 4$ situation is fixed (equal to $\log 2$) in some cases. The behaviour of the Sano-Watatani angle operator is reverse; while it is a singleton set in the $2\times 2$ situation, cardinality of the spectrum of the angle operator is not fixed in the $4\times 4$ situation and it depends on the index of the intersection. Thus, the outcome changes depending on the order of the Hadamard matrices we are dealing with, and whether they are Hadamard equivalent or inequivalent. These contrasts indicate the fact that pairs of Hadamard subfactors are quite interesting and requires extensive investigation. In continuation of \cite{BG1,BG2}, the present paper is the first part of our investigation of pairs of Hadamard subfactors arising from pairs of complex Hadamard matrices in the Hadamard equivalence class of the Fourier matrix. Let us describe it in a slight more precise manner.

Let $F_n$ denote the $n\times n$ Fourier matrix (also called DFT matrix) and consider the Hadamard equivalence class $[F_n]$. If we pick any two matrices $u,v\in[F_n]$, and obtain the Hadamard subfactors $R_u\subset R$ and $R_v\subset R$, by a characterization result in \cite{BG1} we know when can we get $R_u\neq R_v$. This is governed by an equivalence relation, denoted by `$\sim$', among $u\mbox{ and }v$ that is finer than the Hadamard equivalence `$\simeq$'. If $u\nsim v$, then we have $R_u\neq R_v$ and vice versa (Theorem $4.2$, \cite{BG1}). The `Hamming numbers' of the rows of the unitary matrix $u^*v$ play the central role here. Note that in this case, both the Hadamard subfactors are crossed product by outer actions of the cyclic group $\mathbb{Z}_n$ on the hyperfinite type $II_1$ factor $R$. The first natural question that we encounter is whether $R_u\cap R_v$ is a factor, and if so, can we characterize $R_u\cap R_v\subset R$? An astute reader must have noticed that intersection of (finite-index) factors need not be a factor, and even if it is, there is no guarantee that it is of finite-index \cite{Jo1} (also see Section $7.2$ in \cite{BG1}). Indeed, Jones-Xu \cite{JX} showed that finiteness of the Pimsner-Popa index of the intersection is equivalent to the finiteness of the Sano-Watatani angle between the subfactors. While investigating the factoriality of the intersection of two Hadamard sufactors arising from $\{(u,v):u,v\in[F_n];\,u\nsim v\}$, it turns out that, quite interestingly, there is some sort of rigidity present between the cases $n=2,3$ and $n\geq 4$. Since the $n\geq 4$ cases are quite involved and require different techniques than in $n=2,3$, we are sort of compelled to break these two situations. In the present article, which is more or less self-contained, we settle the $n=3$ case (the $n=2$ case is in \cite{BG1}), and the general $n\geq 4$ cases will be dealt separately.

Let us now briefly mention the findings of this article and explain the kind of rigidity that appears. Consider a pair of Hadamard subfactors $R_u\subset R\mbox{ and }R_v\subset R$ arising from complex Hadamard matrices $u\mbox{ and }v$ of order $3$ such that $u\nsim v$ (so that $R_u\neq R_v$). We show that $R_u\cap R_v$ is a subfactor of the hyperfinite type $II_1$ factor $R$. Moreover, quite interestingly, it is a vertex model subfactor of index $9$ in $R$. Thus, the pair $(u,v)$ jointly produces a {\it bi-unitary} permutation matrix (upto equivalence in the sense of \cite{KSV}) of order $9$ through some quantum operation `$u\boxtimes v$'. We compute the relative commutant $(R_u\cap R_v)^\prime\cap R$ and show that $R_v=\mbox{Ad}_w(R_u)$, where $w\in(R_u\cap R_v)^\prime\cap R$ is a unitary. Thus, any pair of Hadamard subfactors of index $3$ are either equal or conjugate by a unitary in the relative commutant (in the hyperfinite type $II_1$ factor $R$) of their intersection. Furthermore, we characterize the subfactor $R_u\cap R_v\subset R$ by identifying it with a vertex model subfactor of the Krishnan-Sunder type \cite{KS}. This subfactor has depth $2$.
\smallskip

The following are the main results of this article.

\begin{thm}
Let $u\mbox{ and }v$ be complex Hadamard matrices of order $3$ such that $u\nsim v$. Consider the corresponding Hadamard subfactors $R_u\subset R\mbox{ and }R_v\subset R$, where $R$ is the hyperfinite type $II_1$ factor. Then, we have the following\,:
\begin{enumerate}[$(i)$]
\item The intersection $R_u\cap R_v$ is a vertex model subfactor of $R$ with $[R:R_u\cap R_v]=9$;
\item $(R_u\cap R_v)^\prime\cap R=\bbc^3$ and $R_u=\mathrm{Ad}_w(R_v)$, where $w\in(R_u\cap R_v)^\prime\cap R$ is a unitary.
\end{enumerate}
\end{thm}

See \Cref{ks type} for the principal graph of $R_u\cap R_v\subset R$.

\begin{thm}
For $u=\mathrm{diag}\{1,e^{i\alpha_1},e^{i\alpha_2}\}F_3$ and $v=\mathrm{diag}\{1,e^{i\beta_1},e^{i\beta_2}\}F_3$, we have the following\,:
\begin{enumerate}[$(i)$]
\item The Pimsner-Popa probabilistic number $\lambda(R_u,R_v)$ is equal to $1/3;$
\item The interior and exterior angle both are equal to $\arccos(|\zeta|^2)$, where $\zeta=\frac{1}{3}\big(e^{i(\alpha_1-\beta_1)}+e^{-i(\alpha_2-\beta_2)}+e^{-i(\alpha_1-\beta_1)}e^{i(\alpha_2-\beta_2)}\big)\in\bbc;$
\item The Sano-Watatani angle between the subfactors $R_u\mbox{ and }R_v$ is the singleton set $\{\arccos|\zeta|\};$
\item $h(R_u|R_v)=\eta\left(\frac{1}{9}\,|1+e^{i(\beta_1-\alpha_1)}+e^{i(\beta_2-\alpha_2)}|^2\right)+\eta\left(\frac{1}{9}\,|1+e^{i(\beta_1-\alpha_1)}\omega+e^{i(\beta_2-\alpha_2)}\omega^2|^2\right)\\
\hspace*{1.9cm} +\eta\left(\frac{1}{9}\,|1+e^{i(\beta_1-\alpha_1)}\omega^2+e^{i(\beta_2-\alpha_2)}\omega|^2\right),$\\
where $\omega$ is a primitive cube root of unity. Furthermore, $h(R_u|R_v)\leq H(R_u|R_v)\leq\log 3$. When the quadruple $(R_u\cap R_v \subset R_u, R_v \subset R)$ is a commuting square, $h(R_u|R_v)= H(R_u|R_v)=-\log\lambda(R_u,R_v)=\log 3$.
\end{enumerate}
\end{thm}

En route, we have characterized when the quadruple $(R_u\cap R_v\subset R_u, R_v\subset R)$ of $II_1$ factors forms a commuting square.

In general for $n\geq 4$, if $u,v\in[F_n]$ with $u\nsim v$, in an upcoming article we shall show that $R_u\cap R_v\subset R$ is always a subfactor but fails to be a vertex model subfactor. Moreover, $(R_u\cap R_v)^\prime\cap R\neq\bbc^n$ when $n\geq 4$. These are in sharp contrasts with the situations $n=2,3$. The present paper establishes the fact that the $n=2$ and $n=3$ cases more or less travel hand in hand, and in an upcoming article we shall explore these rigidity and the quantum operation `$u\boxtimes v$' discussed above.


\newsection{Hadamard subfactors and vertex model subfactors}\label{Sec 2}

\noindent\textbf{Notations:} Throughout the article, we reserve the following notations.
\begin{enumerate}[$(i)$]
\item $M_n$ denotes the algebra of $n\times n$ matrices over $\bbc$. By $\Delta_n$, we denote the diagonal subalgebra (Masa) in $M_n$. Also, $\mathcal{U}(M_n)$ denotes the group of unitary matrices.
\item $\mathrm{diag}\{\mu_1,\ldots,\mu_n\}$ denotes the diagonal matrix in $M_n(\bbc)$ with $\mu_j$'s in the diagonal and zero elsewhere.
\item We often use the shorthand notation $(N\subset P,Q\subset M)$ to denote the following quadruple
\[\begin{matrix}
P &\subset & M\\
\cup & &\cup\\
N &\subset & Q
\end{matrix}\]
of (finite)von Neumann algebras.
\end{enumerate}

Let us start by briefly recalling the construction of Hadamard subfactors, also called spin model subfactors, and the vertex model subfactors from \cite{JS}. These are obtained from complex Hadamard matrices.
\begin{dfn}
A complex Hadamard matrix $H$ is a $n\times n$ matrix with complex entries of the same modulus such that $HH^{*}=nI_{n}$.
\end{dfn}

Notice that $\frac{1}{\sqrt{n}}H$ is a unitary matrix. In this article, we work with the following definition of complex Hadamard matrices, as is customary in the world of subfactors.
\begin{dfn}
A complex Hadamard matrix of order $n$ is a unitary matrix such that each of its entry has the modulus $1/\sqrt{n}$.
\end{dfn}
We denote by $F_n$ the Fourier matrix (also called DFT matrix) $\big(\omega^{jk}/\sqrt{n}\big)_{j,k=0,\ldots,n-1}$, where $\omega=e^{-2\pi i/n}$ is a primitive $n$-th root of unity. Two complex Hadamard matrices are Hadamard equivalent, to be denoted by $H_{1}\simeq H_{2}$, if there exist diagonal unitary matrices $D_1, D_2$ and permutation matrices $P_1, P_2$ such that 
\begin{eqnarray}\label{Hequivalence}
H_{1}=D_{1}P_{1}H_{2}P_{2}D_{2}.
\end{eqnarray}
It is known that for $n=2,3,5$, all complex Hadamard matrices are Hadamard equivalent to the Fourier matrix $F_n$. However, complete classification of complex Hadamard matrices is not known and quite hard beyond $n=5$.

Let $u$ be a complex Hadamard matrix of order $n$. It is known \cite{GHJ,JS} that the following quadruple
\[
\begin{matrix}
\mbox{Ad}_u(\Delta_n) &\subset & M_n\\
\cup & &\cup\\
\bbc &\subset & \Delta_n
\end{matrix}
\]
is a non-degenerate (also called symmetric) commuting square. Iterating Jones' basic construction, we obtain the spin model subfactor $R_u\subset R$ of the hyperfinite type $II_1$ factor $R\,$: 
\[
\begin{matrix}
\Delta_n &\subset & M_n  &\subset^{\,e_1} & \Delta_n\otimes M_n  &\subset^{\,e_2} &\ldots \ldots &\subset  R\\
\cup & & \cup  & &\cup &   &  &\quad\cup\\
\bbc &\subset &\mbox{Ad}_u(\Delta_n)  &\subset &\mbox{Ad}_{u_1}(M_n) &\subset &\ldots \ldots  &\subset  R_u\\
\end{matrix}
\]
where $e_j$'s are the Jones' projections for the basic construction of $\Delta_n\subset M_n$, and $u_j$'s are certain unitary matrices given by the following :
\begin{ppsn}[\cite{N}]\label{unitaries}
Let $u=(u_{ij})$ and $D_u=\displaystyle{\sqrt{n}\sum_{i=1}^{n}\sum_{j=1}^{n}\overline{u_{ij}}(E_{ii}\otimes E_{jj})}$. Then,
\begin{center}
$u_{2k+1}=(I_n \otimes u_{2k})(D_u \otimes I^{(k)}_n)\quad\mbox{and}\quad u_{2k}=u_{2k-1}(u \otimes I_{n}^{(k)})$
\end{center}
are the unitary matrices in the tower of the basic construction for the Hadamard subfactor $R_u\subset R$.
\end{ppsn}

At first glance, it may appear that the unitary matrices described above differs from that in \cite{N}. However, this is only due to our choice of embedding. In \cite{N}, the embedding $M_n\subset M_n\otimes M_n$ is on the right, that is, $x\mapsto x\otimes I_n$ (i,e., $M_n\otimes \bbc\subset M_n\otimes M_n$); whereas in our convention it is $x\mapsto I_n\otimes x$, that is, $\bbc\otimes M_n\subset M_n\otimes M_n$ with the embedding $x\mapsto\mbox{bl-diag}\{x,\ldots,x\}$. Throughout the article, we adhere to this convention.

Not much is known about this class of subfactors; however, it is known that these are irreducible with second relative commutant abelian. Now, we discuss the vertex model subfactors.
\begin{dfn}[\cite{JS}]
A unitary matrix $u=(u^{\alpha a}_{\beta b})$ in $M_n\otimes M_k$ is said to be a bi-unitary matrix if the block-transpose $\widetilde{u}=(\widetilde{u}^{\alpha a}_{\beta b})$, defined by $\widetilde{u}^{\alpha a}_{\beta b}:=u^{\beta a}_{\alpha b}$, is also a unitary matrix in $M_n\otimes M_k$.
\end{dfn}
It is known that the following quadruple
\[
\begin{matrix}
\mbox{Ad}_u(M_n\otimes\bbc) &\subset & M_n\otimes M_k\\
\cup & &\cup\\
\bbc &\subset & \bbc\otimes M_k
\end{matrix}
\]
is a non-degenerate commuting square precisely when $u$ is a bi-unitary matrix \cite{JS}. Iterating Jones basic construction, we obtain a subfactor $R_u\subset R$ of the hyperfinite $II_1$ factor $R$ such that $[R:R_u]=k^2$. Unlike the spin model subfactors, the vertex model subfactors need not be irreducible. There is a natural equivalence relation on the set of bi-unitary matrices in $M_n\otimes M_k$ such that the subfactors arising from equivalent bi-unitary matrices are conjugate. This equivalence relation is given by $u_1\sim u_2$ if and only if there exist unitary matrices $a,c\in M_n$ and $b,d\in M_k$ such that $u_1=(a\otimes b)u_2(c\otimes d)$ (see Section $4$, \cite{KSV} and Section $2$, \cite{KS}).

In this article, we only need vertex model subfactors arising from permutation bi-unitary matrices investigated in \cite{KSV, KS}. Let us briefly recall few essential facts needed in this article.

\begin{lmma}[Lemma $1$, \cite{KS}]\label{bi-unitary}
Let $\Omega_n=\{1,2,3,\ldots,n\}$ and $u\in M_n\otimes M_n$. The following are equivalent.
\begin{enumerate}[$(i)$]
\item u is bi-unitary permutation matrix;
\item there exist permutations $\{\rho_k\, :\, k \in \Omega_n\}\subseteq S(\Omega_n)$ and $\{\lambda_k\, :\, k \in \Omega_n \}\subseteq S(\Omega_n)$, where $S(\Omega_n)$ denotes the group of all permutations of $\Omega_n$, such that
\begin{enumerate}[$(a)$]
\item the equation
\[
\pi(j,\ell)=(\rho_\ell(j),\, \lambda_j(\ell))
\]
defines a permutation $\pi \in S(\Omega_n \times \Omega_n)$; and
\item $u^{ik}_{j\ell}=\delta_{(i,k),\pi(j,\ell)}=\delta_{i,\rho_\ell(j)}\delta_{k,\lambda_j(\ell)}.$
\end{enumerate}
\end{enumerate}
\end{lmma}

\begin{dfn}[Definition $2$, \cite{KS}]\label{Definition}
Define
\begin{center}
$P_n:=\{\pi\in S(\Omega_n\times\Omega_n)\,:\exists\,\,\lambda,\rho:\Omega_n\rightarrow S(\Omega_n)\,\,\mbox{ such that }\,\pi(j,\ell)=(\rho_\ell(j),\,\lambda_j(\ell))$
\end{center}
\begin{center}
$\hspace*{10.5cm}\mbox{ for all }\ell,j\in\Omega_n\}\,,$
\end{center}
where $\lambda_j$ (resp., $\rho_\ell$) denotes the image of $j$ (resp., $\ell$) under the map $\lambda$ (resp., $\rho$).
\end{dfn}

\Cref{bi-unitary} and \Cref{Definition} show that there exists a bijection between bi-unitary permutation matrices of size $n^2$ and elements $\pi\leftrightarrow (\lambda, \rho) \in P_n$,  given by $u^{ik}_{j\ell}=\delta_{(i,k),\pi(j,\ell)}=\delta_{i, \rho_\ell(j)}\delta_{k, \lambda_j(\ell)}.$

In particular, for $n=3$ there exist $18$ inequivalent bi-unitary permutation matrices. Among them, the principal graphs corresponding to $17$ bi-unitary permutation matrices are obtained in \cite{KS}, and the remaining one is obtained in \cite{BINA}. 


\newsection{Pairs of Hadamard subfactors}\label{Sec 3}

Let $u\mbox{ and }v$ be two distinct complex Hadamard matrices of order $n$. By the construction of Hadamard subfactors discussed in \Cref{Sec 2}, we obtain $R_u\subset R$ and $R_v\subset R$. Although $u\neq v$, it may very well happen that $R_u=R_v$. In \cite{BG1}, $R_u\neq R_v$ is completely characterized. To achieve this, for complex Hadamard matrices $u$ and $v$ in $M_n(\bbc)$, define $u\sim v$ if there exists a permutation matrix $P\in M_n$ and a diagonal unitary matrix $D\in M_n$ such that $v=uPD$. Then, the equivalence relation `$\sim$' is finer than the Hadamard equivalence relation defined in \Cref{Hequivalence}. We have the following characterization result.

\begin{thm}[Theorem $4.2$, \cite{BG1}]\label{subrelation}
\begin{enumerate}[$(i)$]
\item For distinct $n\times n$ complex Hadamard matrices $u\mbox{ and }v$, the pair of Hadamard subfactors $R_u\subset R\mbox{ and }R_v\subset R$ are distinct $(\mbox{that is, }R_u\neq R_v)$ if and only if $\,u\nsim v$.
\item If two $n\times n$ complex Hadamard matrices $u\mbox{ and }v$ are Hadamard inequivalent, then the corresponding spin model subfactors $R_u\mbox{ and }R_v$ of $R$ are always distinct $(R_u\neq R_v)$.
\end{enumerate}
\end{thm}

\begin{figure}
    \centering
    \begin{tikzpicture}
        \draw (0,0) node [scale =0.85]{{${\Delta}_n$}} (-4,-1) node [scale =0.85]{{$\bbc$}} (4,-1) node [scale =0.85]{{$\bbc$}};
        \draw[->-=1] (-3.8, -.9)--(-1,-.15);
        \draw[->-=1] (3.8,-.9)--(1,-.15);
        \draw[->-=1](0,.5)--(0,1.5); 
      
        \draw[->-=1](-4,-1+.5-.2)--(-4,-1+2-.5);
        \draw[->-=1](4,-1+.5-.2)--(4,-1+2-.5);
        \draw (0,0+2) node [scale =0.85] {{$M_n$}} (-4,-1+2) node [scale =0.85]{{$u\Delta_nu^*$}} (4,-1+2) node [scale =0.85] {{$v\Delta_nv^*$}};
        \draw[->-=1] (-3.5,-.8+2-.2)--(-1,-.15+2);
        \draw[->-=1] (3.5,-.8+2-.2)--(1,-.15+2);
        \draw[->-=1](0,.5+2)--(0,1.5+2); 
      
        \draw[->-=1](-4,-1+.5+2)--(-4,-1+2-.5+2);
        \draw[->-=1](4,-1+.5+2)--(4,-1+2-.5+2);
        \draw (0,0+4) node[scale =0.85] {{$\Delta_n\otimes M_n$}} (-4,-1+2+2) node [scale =0.85] {{$u_{1} M_n u^*_{1}$}} (4,-1+2+2) node [scale =0.85]{{$v_{1}M_nv^*_{1}$}}; 
       \draw[->-=1] (-3.4,-.8+2+2-.2)--(-1,-.15+2+2);
        \draw[->-=1] (3.4,-.8+2+2-.2)--(1,-.15+2+2);
        \draw[->-=1](0,.5+2+2)--(0,1.5+2+2);    
         \draw[->-=1] (-2.9,-.9+2+2+2)--(-1,-.15+2+2+2);
        \draw[->-=1] (2.9,-.8+2+2+2-0.1)--(1,-.15+2+2+2);
         \draw[->-=1](-4,-1+.5+2+2)--(-4,-1+2-.5+2+2);
        \draw[->-=1](4,-1+.5+2+2)--(4,-1+2-.5+2+2);
         \draw (0,0+6) node[scale =0.80] {{$M_n\otimes M_n$}} (-4,-1+2+2+2) node [scale=0.80]{{$u_{2}(\Delta_n\otimes M_n) u^*_{2}$}} (4,-1+2+2+2) node [scale =0.80] {{${v}_{2}(\Delta_n\otimes M_n) v^*_{2}$}};
        \draw[ dashed](0,.5+2+4)--(0,1.5+2+4); 
        \draw[dashed](-4,-1+.5+2+4)--(-4,-1+2-.5+2+4);
        \draw[dashed](4,-1+.5+2+4)--(4,-1+2-.5+2+4);
        \draw (0,0+7.9) node [scale =0.85] {{$R$}} (-4,-1+2+6-.1) node [scale =0.85]{{$R_u$}} (4,-1+2+6-.1) node [scale =0.85]{{$R_v$}};  
          \draw[->-=1] (-3.4,7)--(-.8,-.15+8-.1);
        \draw[->-=1] (3.4,7)--(.8,-.15+8-.1);
    \end{tikzpicture}
    \caption{Construction of pairs of Hadamard subfactors}\label{main fig}
\end{figure}
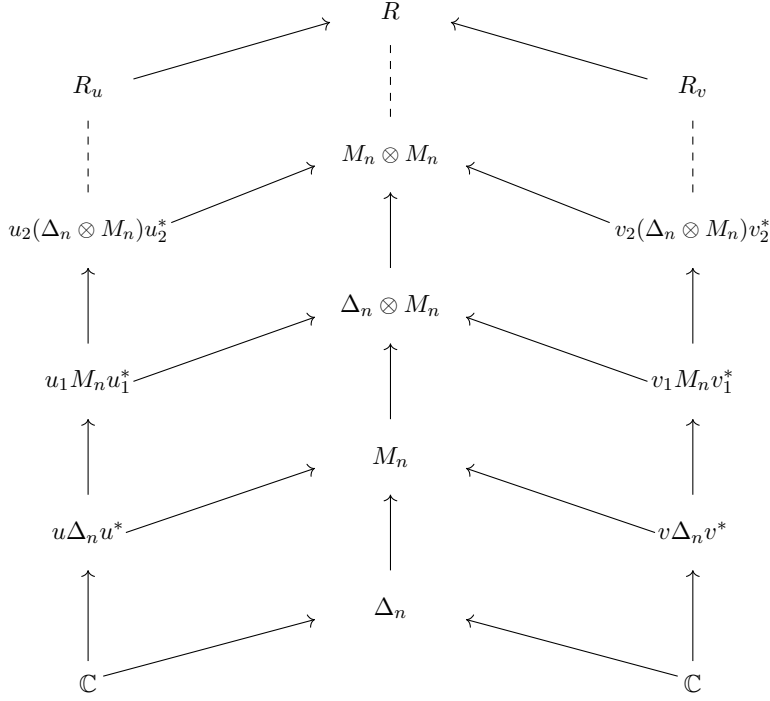

The construction of pair of Hadamard subfactors is depicted in \Cref{main fig} (see \cite{BG1} for detail). Note that `$B\subset A$' has been marked by `$B\rightarrow A$' in the figure. The first difficulty in investigating two Hadamard subfactors lies in the fact that it is not automatic to make the diagram $R_u\subset R\supset R_v$ a quadruple of factors. More precisely, it is not guaranteed whether there exists a factor $N$ such that $N\subseteq R_u\cap R_v$.  Moreover, even if $R_u\cap R_v$ becomes a factor, it may be of infinite index, in which case the quadruple still remains a challenge (Theorem $7.30$, \cite{BG1}).

\begin{figure}
  \centering
  \begin{tikzpicture}
        \draw[thick,-<-=.5] (0,0)--(4,0);
        \draw[thick,->-=.7] (4,0)--(4,2);
        \draw[thick,->-=.5] (4,2)--(0,2);
        \draw[thick,-<-=.5] (0,2)--(0,0);
        \draw[thick,-<-=.5] (1,1)--(5,1);
        \draw[thick,->-=.5] (5,1)--(5,3);
        \draw[thick,->-=.5] (5,3)--(1,3);
        \draw[thick,-<-=.7] (1,3)--(1,1);
        \draw[thick,->-=.5] (0,0)--(1,1);
        \draw[thick,->-=.5] (4,0)--(5,1);
        \draw[->-=.5] (4,0)--(1,1);
        \draw[thick,->-=.5] (4,2)--(5,3);
        \draw[thick,->-=.5] (0,2)--(1,3);
        \draw[->-=.5] (4,2)--(1,3);

 \draw (-.3,-.2) node {$B^1_0$};
        \draw (-.3,2.1) node {$B^1_1$};
        \draw (5.2, 1) node {$B^2_0$};
         \draw (.7, 1.2) node{$A_0$};
        \draw (.7,3.2) node{$A_1$};
        \draw (5.25, 2.9) node {$B^2_1$};
        \draw (4.4,1.9) node {$C_1$};
        \draw (4.25, -.1) node {$C_0$};
\end{tikzpicture}
\caption{commuting cube}\label{com1}
\end{figure}
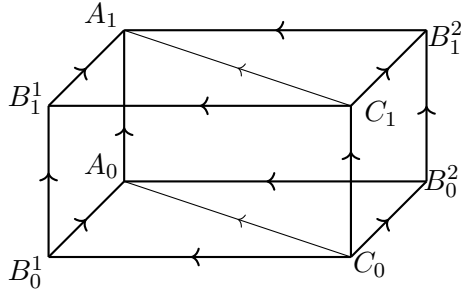

A natural and relevant concept to investigate pairs of subfactors, or `two subfactors', is the notion of ``commuting cube" introduced in Section $3$, \cite{BG1} (similar notion has appeared in different context in the finite-dimensional situation \cite{K}) which may be thought of as one-dimension higher object than commuting square. This concept has been extensively used to deal with the question of factoriality of $R_u\cap R_v$ investigated in \cite{BG1}. We recall it briefly here as we shall need it in this article. Consider the cube of finite von Neumann algebras depicted in \Cref{com1}, where $A_1$ is equipped with a faithful normal tracial state, and $C_0=B_0^1\cap B_0^2,\,C_1=B_1^1\cap B_1^2$. In \Cref{com1}, `$B\rightarrow A$' denotes $B\subset A$.

\begin{ppsn}[\cite{BG1}]\label{slice}
Suppose that the adjacent faces $(B_0^j\subset A_0,B_1^j\subset A_1),\,j=1,2,$ in \Cref{com1} are commuting squares. Then, $C_1\cap B_0^j=C_0$ for $j=1,2,$ and $C_1\cap A_0=C_0$. Moreover, the remaining faces $(C_0\subset B_0^j\,,\,C_1\subset B_1^j),\,j=1,2$, and the slice $(C_0\subset A_0\,,\,C_1\subset A_1)$ are also commuting squares.
\end{ppsn}

\begin{dfn}[\cite{BG1}]\label{commcube}
A cube of finite von Neumann algebras as in \Cref{com1} is called a commuting cube if both the adjacent faces $(B_0^j\subset A_0,B_1^j\subset A_1)$ for $j=1,2$ are a commuting square.
\end{dfn}

\begin{rmrk}\rm
\begin{enumerate}[$(i)$]
\item The floor $(C_0\subset B_0^1,B_0^2\subset A_0)$ and the roof $(C_1\subset B_1^1,B_1^2\subset A_1)$ in a commuting cube need not be a commuting square.
\item If the adjacent faces in a commuting cube are non-degenerate commuting squares, then although the slice $(C_0\subset A_0,C_1\subset A_1)$ is a commuting square, it may fail to become non-degenerate.
\end{enumerate}
\end{rmrk}

\begin{ppsn}[\cite{BG1}]\label{klz}
Suppose that we have a commuting cube as in \Cref{com1}. If the roof $(C_1\subset B_1^1,B_1^2\subset A_1)$ is a commuting square, then the floor $(C_0\subset B_0^1,B_0^2\subset A_0)$ is also a commuting square. The converse need not hold.
\end{ppsn}

\begin{dfn}[\cite{BG1}]\label{nondegeneratecommcube}
A commuting cube in \Cref{com1} is called a non-degenerate commuting cube if both the adjacent faces and the slice are non-degenerate commuting squares.
\end{dfn}

Recall the basic construction of non-degenerate commuting cube discussed in Section $3$, \cite{BG1}. The following theorem will be needed in this article to compute the Sano-Watatani angle operator.

\begin{thm}[\cite{BG1}]\label{gggg}
Suppose that $(N\subset P,Q\subset M)$ is a quadruple of $II_1$ factors obtained as an iterated basic construction of a non-degenerate commuting cube of finite-dimensional algebras depicted in \Cref{com1}. Consider the following nonnegative matrix
\[
S_0:=E^{A_0}_{B^2_0}E^{A_0}_{B^1_0}E^{A_0}_{B^2_0}-E^{A_0}_{C_0}
\]
and suppose that $S_0\neq 0$ with $S_0^{\,2}=\alpha S_0$ for some $\alpha\in\mathbb{R}_+$. Then, $(N\subset P,Q\subset M)$ is not a commuting square, and moreover $\mathrm{Ang}_M(P,Q)$ is the singleton set $\{\mathrm{arccos}\sqrt{\alpha}\}$.
\end{thm}


\newsection{Hadamard subfactors of index \texorpdfstring{$3$}{} and characterization of the intersection}\label{Sec 4}

In this section, we focus on complex Hadamard matrices of order $3$ and consider the pair of Hadamard subfactors each of index $3$. Our goal is to prove the factoriality of the intersection, compute the relative commutant, and characterize the intersection. En route, we prove the conjugacy of the Hadamard subfactors.

It is known that there is unique (Hadamard) equivalence class in this case, and any $3\times 3$ complex Hadamard matrix is Hadamard equivalent to the Fourier matrix $F_3$. If $(u,v)$ is any pair of such matrices, and $R_u\subset R,\,R_v\subset R$ are the corresponding Hadamard subfactors, then $R_u=R_v$ if and only if $u\sim v$ by \Cref{subrelation}. Therefore, in order to obtain a pair of Hadamard subfactors, we must choose $u,v$ satisfying $u\nsim v$, and obtain the following :
\begin{eqnarray}\label{pxz}
\begin{matrix}
R_u &\subset & R\\
 & & \cup\\
& & R_v
\end{matrix}
\end{eqnarray}
through the construction depicted in \Cref{main fig}. Note that at present \Cref{pxz} can be made a quadruple of von Neumann algebras only, since the factoriality of $R_u\cap R_v$ is not clear. Employing the definition of Hadamard equivalence, if we write $u=D_1P_1F_3P_2D_2$ and $v=\widetilde{D}_1\widetilde{P}_1F_3\widetilde{P}_2\widetilde{D}_2$ (see \Cref{Hequivalence}), then due to \Cref{subrelation} we see that $R_u=R_w$, where $w=D_1P_1F_3$ (as $u=wP_2D_2$). Therefore, without loss of generality, we can assume that $u=D_1P_1F_3$, and similarly $v=\widetilde{D}_1\widetilde{P}_1F_3$. The following result shows that we can say even more. Fix the following diagonal matrices throughout the rest of the paper.
\begin{eqnarray}\label{pq}
\mathscr{D}_1=I_3=\begin{bmatrix}
    1 & 0 & 0 \\
    0 & 1 & 0 \\
    0 & 0 & 1 
\end{bmatrix}\quad,\quad
\mathscr{D}_2=\begin{bmatrix}
    1 & 0 & 0 \\
    0  & \omega & 0\\
    0 & 0& \omega^2
\end{bmatrix}\quad,\quad
\mathscr{D}_3=\begin{bmatrix}
    1 & 0 & 0 \\
    0  & \omega^2 & 0\\
    0 & 0& \omega
\end{bmatrix}\,\,.
\end{eqnarray}

\begin{thm}\label{equivalence}
For any $3\times 3$ complex Hadamard matrix $u$, there exists a diagonal unitary matrix $D$ such that $u\sim D F_3$. Furthermore, $D_1F_3\sim D_2F_3$ if and only if $D^*_2D_1$ is equal to scalar multiple (of modulus one) of any of the diagonal matrices $\mathscr{D}_j,\,j=1,2,3,$ defined in \Cref{pq}.
\end{thm}
\begin{prf}
Let $S_3$ be the subgroup of $GL_3(\bbc)$ consisting of all $3\times 3$ permutation matrices. Consider the following subset of $S_3$
\[
G=\{P\in S_3\,:\,PF_3\sim F_3\}\,.
\]
It is easy to check that $G$ is a subgroup of $S_3$. We claim that $G=S_3$. Since $S_3$ can be generated by a $2$ cycle and a $3$ cycle, we take $(23)$ and $(132)$ as generators. Therefore, the matrices $\gamma_1=E_{11}+E_{23}+E_{32}$ and $\gamma_2=E_{12}+E_{23}+E_{31}$ generate $S_3$ in $GL_3(\bbc)$. To fulfil the claim, it is enough to show that $\gamma_1,\gamma_2\in G$. Now, it is a straightforward verification that $\gamma_1F_3=F_3\gamma_1$ and $\gamma_2F_3=F_3\,\mbox{diag}\{1,\omega^2,\omega\}$. Therefore, $\gamma_1F_3\sim F_3$ and $\gamma_2F_3\sim F_3$, which completes the proof of the claim.

Now, for any complex Hadamard matrix $u\in M_3$, write $u=D_1P_1F_3P_2D_2$, where $D_1,D_2$ are diagonal unitary matrices and $P_1,P_2\in S_3$. Hence, we have $u\sim D_1P_1F_3$. Since $G=S_3$, we have $P_1F_3 \sim F_3$. Therefore, we have $u\sim D_1F_3$, which completes the proof of the first part.
\smallskip

For the second part, suppose that $D_2^*D_1\in\{\alpha\mathscr{D}_j:\alpha\in\mathbb{S}^1,\,j=1,2,3\}$. Consider the alternating subgroup $A_3=\{\mbox{id},(123),(132)\}\unlhd S_3$ of even permutations. The matrix representation of $A_3$ is $\{I_3,\sigma_1=(123)=E_{21}+E_{32}+E_{13},\sigma_2=(132)=E_{31}+E_{12}+E_{23}\}$. Observe that $\mathscr{D}_1F_3=F_3,\,\mathscr{D}_2F_3=F_3\sigma_1$ and $\mathscr{D}_3F_3=F_3\sigma_2$. That is, $\{\mathscr{D}_jF_3:j=1,2,3\}=\{F_3,F_3\sigma_1,F_3\sigma_2\}$. Now, if $D_2^*D_1=\alpha\mathscr{D}_j$ for some $j\in\{1,2,3\}$, then we get $D_1F_3=\alpha D_2\mathscr{D}_jF_3=D_2F_3P(\alpha I_3)$ for some permutation matrix $P\in A_3=\{I_3,\sigma_1,\sigma_2\}$. Therefore, we have $D_1F_3\sim D_2F_3$. Conversely, suppose that $D_1F_3\sim D_2F_3$. Then, $D_1F_3=D_2F_3PD$, where $P$ is a permutation matrix and $D$ is a diagonal unitary matrix. Then, $D_2^*D_1F_3=F_3PD$. Since $F_3$ is the Fourier matrix, for any permutation $P\in S_3$, the first column of $F_3P$ can be $(1\,\,1\,\,1)$ or $(1\,\,\omega\,\,\omega^2)$ or $(1\,\,\omega^2\,\,\omega)$. Writing $D=\mbox{diag}\{z_1,z_2,z_3\}$, where $z_j\in\mathbb{S}^1$, and $D_2^*D_1=\mbox{diag}\{\lambda_1,\lambda_2,\lambda_3\}$, we see that only the following combinations are possible
\begin{center}
$(\lambda_1,\lambda_2,\lambda_3)=z_1(1,1,1)\,\,,\,\,(\lambda_1,\lambda_2,\lambda_3)=z_1(1,\omega,\omega^2)\,\,,\,\,(\lambda_1,\lambda_2,\lambda_3)=z_1(1,\omega^2,\omega)\,.$
\end{center}
Since $z_1\in\mathbb{S}^1$, we have $D_2^*D_1\in\{\alpha\mathscr{D}_j:\alpha\in\mathbb{S}^1,\,j=1,2,3\}$.\qed
\end{prf}

As explained at the begining of this section, without loss of generality, we can take $u=D_1P_1F_3$ and $v=D_2P_2F_3$. Therefore, due to \Cref{equivalence}, we can further discard  permutation matrices $P_1\mbox{ and }P_2$ (since $R_u=R_{DF_3}$ by \Cref{subrelation}), and work with the following set-up.
\medskip

\noindent\textbf{The set-up:}~Throughout the rest of the paper, we take $u=D_1F_3\mbox{ and }v=D_2F_3$ such that $D_2^*D_1\notin\{\alpha\mathscr{D}_j:\alpha\in\mathbb{S}^1,\,j=1,2,3\}$ (see \Cref{pq} for the notations $\mathscr{D}_j$).
\smallskip

\begin{notation}\label{not}\rm
For $k \in \mathbb{N} \cup \{0\}$, define
\begin{align*}
        A_{2k}:=M_3\otimes M^{(k)}_3\,,\quad& A_{2k+1}:=\Delta_3\otimes M_3 \otimes M^{(k)}_3\,,\\
        B^{u}_{2k}:=\mbox{Ad}_{u_{2k}}\big(\Delta_3\otimes M_3^{(k)}\big)\,,\quad& B^{u}_{2k+1}:=\mbox{Ad}_{u_{2k+1}}\big(M_3 \otimes M_3^{(k)}\big)\,,\\
        B^{v}_{2k}:=\mbox{Ad}_{v_{2k}}\big(\Delta_3\otimes M_3^{(k)}\big)\,,\quad& B^{v}_{2k+1}:=\mbox{Ad}_{v_{2k+1}}\big(M_3 \otimes M_3^{(k)}\big)\,,\\
        C_{2k}:=B^u_{2k}\cap B^v_{2k}\,,\quad& C_{2k+1}:=B^u_{2k+1}\cap B^v_{2k+1}.
\end{align*}
These are the even, respectively odd, steps in the tower of basic construction depicted in \Cref{main fig}. The above notations will be reserved throughout the article.
\end{notation}
Note that $R=\overline{\cup A_{2k}}^{\,\mbox{sot}}=\overline{\cup A_{2k+1}}^{\,\mbox{sot}},\,R_u=\overline{\cup B^u_{2k}}^{\,\mbox{sot}}=\overline{\cup B^u_{2k+1}}^{\,\mbox{sot}},\,R_v=\overline{\cup B^v_{2k}}^{\,\mbox{sot}}=\overline{\cup B^v_{2k+1}}^{\,\mbox{sot}}$ and $R_u\cap R_v=\overline{\cup\, C_{2k}}^{\,\mbox{sot}}=\overline{\cup\, C_{2k+1}}^{\,\mbox{sot}}$ (see Sections $6\mbox{ and }7$, \cite{BG1} for detail on this).


\subsection{Factoriality of \texorpdfstring{$R_u\cap R_v$}{}}

Goal of this subsection is to prove factoriality of $R_u\cap R_v$. We start by fixing notations for some specific matrices that are used throughout this subsection and later.
\begin{notation}\label{notation}
\begin{enumerate}[$(i)$]
\item Let $\sigma_1:=E_{21}+E_{32}+E_{13}$ and $\sigma_2:=E_{31}+E_{12}+E_{23}$. These matrices are the representations of the alternating subgroup $A_3=\{\mbox{id},(123),(132)\}\unlhd S_3$ of even permutations.
\item For any $y\in M_3$, let $L_y$ denote the left multiplication operator acting on $M_3$, that is, $L_y(x)=yx$ for $x\in M_3$. Define the following three operators acting on $M_3$
\[
Q_0=E^{M_3}_{\Delta_3}\,,\,\,Q_1=L_{\sigma_2}E^{M_3}_{\Delta_3}L_{\sigma_1}\,,\,\,Q_2=L_{\sigma_1}E^{M_3}_{\Delta_3}L_{\sigma_2},
\]
where $E^{M_3}_{\Delta_3}$ is the unique trace preserving conditional expectation onto $\Delta_3$. Thus, for each $x=(x_{ij})_{1 \leq i, j \geq 3}$ in $M_3$, we have the following generalized permutation matrices
\[
Q_0(x)=\begin{bmatrix}
    x_{11} & 0 & 0 \\
    0  & x_{22} & 0\\
    0 & 0&  x_{33}
    \end{bmatrix}\,,\quad Q_1(x)=\begin{bmatrix}
    0 & x_{12} & 0 \\
    0  & 0 & x_{23}\\
    x_{31} & 0&  0
\end{bmatrix}\,,\quad Q_2(x)=\begin{bmatrix}
    0 & 0 & x_{13} \\
    x_{21}  & 0 & 0\\
    0 & x_{32}&  0
\end{bmatrix}.
\]
\item Let $W_2=\mbox{bl-diag}\{I_3,\sigma_1,\sigma_2\}$ and for $k\geq 2$, define $W_{2k}=\displaystyle{\prod_{n=0}^{k-1}I^{(k-1-n)}_{3}\otimes W_2\otimes I^{(n)}_3} \in \Delta_3 \otimes  M^{(k)}_3$. Note that each $W_{2k}$ is a unitary matrix.
\end{enumerate}
\end{notation}

\begin{lmma}\label{Intersection between von Neumann algebras}
For any $k \in \mathbb{N}$, we have the following
\[
\mathrm{Ad}_{u_{2k}^*v_{2k}}\big(\Delta_3 \otimes M^{(k)}_3\big)=\mathrm{Ad}_{W_{2k}}\big(\mathrm{Ad}_{u^*v}(\Delta_3)\otimes M^{(k)}_3\big),
\]
where $W_{2k}$'s are as defined in \Cref{notation}.
\end{lmma}
\begin{prf}
Recall the tower of basic construction from \Cref{unitaries}, along with the corresponding notations there. We claim that $u^*_{2k}v_{2k}=\mathrm{Ad}_{W_{2k}}\big(u^*v \otimes I_3^{(k)}\big)$ for any $k\in\bbn$. For $k=1$, first observe the following
\begin{eqnarray*}
u^*_2v_2 &=& ((I_3\otimes u)D_u(u\otimes I_3))^*(I_3\otimes v)D_v(v\otimes I_3)\nonumber\\ 
&=& (u^* \otimes I_3)D^*_{u}(I_3\otimes u^*v)D_v(v\otimes I_3)\nonumber\\
&=& (F^*_3\otimes I_3)D^*_{{F}_{3}}(I_3\otimes u^*v)D_{F_3}(F_3\otimes I_3)\,.
\end{eqnarray*}
Now, for any $x=(x_{ij}) \in M_3$, it is a straightforward verification that
\[
(F^*_3\otimes I_3)D^*_{{F}_{3}}(I_3\otimes x)D_{F_3}(F_3\otimes I_3)=\begin{bmatrix}
    Q_0(x) & Q_1(x) &  Q_2(x) \\
    Q_2(x) &  Q_0(x) & Q_1(x)\\
    Q_1(x) &  Q_2(x)&  Q_0(x)
\end{bmatrix}
\]
(see \Cref{notation} for the $Q_j(x)$'s). Hence, for $x=u^*v$ we have the following
\[u^*_2v_2=(F^*_3\otimes I_3)D^*_{{F}_{3}}(I_3\otimes u^*v)D_{F_3}(F_3\otimes I_3)=W_2(u^*v\otimes I_3)W^*_2,
\]
which is the basis step of the induction. Assume that the claim is true up to the $k$-th step for some $k\in\bbn$. Then, for the $(k+1)$-th step, using the induction hypothesis we have the following
\begin{eqnarray*}
& & u^*_{2k+2}v_{2k+2}\nonumber\\
&=& \big((I_3 \otimes u_{2k})(D_u \otimes I^{(k)}_3)(u \otimes I^{(k+1)}_3)\big)^*\big((I_3 \otimes v_{2k})(D_v \otimes I^{(k)}_3)(v \otimes I^{(k+1)}_3)\big)\nonumber\\
&=& (u^* \otimes I^{(k+1)}_3)(D^*_u \otimes I^{(k)}_3)(I_3 \otimes u^*_{2k}v_{2k})(D_v \otimes I^{(k)}_3)(v \otimes I^{(k+1)}_3)\quad\qquad(\mbox{ by induction hypothesis})\nonumber\\
&=& (u^* \otimes I^{(k+1)}_3)(D^*_u \otimes I^{(k)}_3)\mathrm{Ad}_{(I_{3} \otimes W_{2k})}\big(I_3 \otimes u^*v\otimes I^{(k)}_3\big)(D_v \otimes I^{(k)}_3)(v \otimes I^{(k+1)}_3)\nonumber\\
&=& \mathrm{Ad}_{(I_3 \otimes W_{2k})}\big((u^* \otimes I^{(k+1)}_3)(D^*_u \otimes I^{(k)}_3)(I_3 \otimes u^*v\otimes I^{(k)}_3)(D_v \otimes I^{(k)}_3)(v \otimes I^{(k+1)}_3)\big)\nonumber\\
&=& \mathrm{Ad}_{(I_3 \otimes W_{2k})}\big(u^*_2v_2 \otimes I^{(k)}_3\big)\nonumber\\
&=& \mathrm{Ad}_{(I_3 \otimes W_{2k})}\big((W_2 \otimes I^{(k)}_3)(u^*v\otimes I^{(k+1)}_3)(W^*_2 \otimes I^{(k)}_3)\big)\nonumber\\
&=& \mathrm{Ad}_{(I_3 \otimes W_{2k})}\mathrm{Ad}_{(W_2 \otimes I^{(k)}_3)}\big(u^*v\otimes I^{(k+1)}_3\big)\nonumber\\
&=& \mathrm{Ad}_{W_{2k+2}}\big(u^*v\otimes I^{(k+1)}_3\big)\,.
\end{eqnarray*}
Since $W_{2k}\in\Delta_3\otimes M^{(k)}_3$ for any $k\in\bbn$, the result follows.\qed
\end{prf}

\begin{lmma}\label{intersection}
For any $k\in\mathbb{N}$, we have the following
\begin{eqnarray*} 
\mathrm{Ad}_{u_{2k}}\big(\Delta_3\otimes M^{(k)}_3\big)\,\bigcap\,\mathrm{Ad}_{v_{2k}}\big(\Delta_3\otimes M^{(k)}_3\big) &=& \mathrm{Ad}_{u_{2k}W_{2k}}\big(\big(\Delta_3\cap\mathrm{Ad}_{u^*v}(\Delta_3)\big)\otimes M^{(k)}_3\big)\nonumber\\
&=&\mathrm{Ad}_{u_{2k}W_{2k}}(\mathbb{C}\otimes M^{(k)}_3),
\end{eqnarray*}
where $W_{2k}$'s are as defined in \Cref{notation}.
\end{lmma}
\begin{prf}
The first equality follows from \Cref{Intersection between von Neumann algebras}. Now, recall that $u=D_1F_3\mbox{ and }v=D_2F_3$, where $D_1\mbox{ and }D_2$ satisfy the condition in \Cref{equivalence}. For the second equality, first observe that
\[
\Delta_3\cap\mathrm{Ad}_{u^*v}(\Delta_3)=\mathrm{Ad}_{F_3^*}\big(\mathrm{Ad}_{F_3}(\Delta_3)\cap\mathrm{Ad}_{D_1^*D_2F_3}(\Delta_3)\big)\,.
\]
Since $D_2^*D_1\notin\{\alpha\mathscr{D}_j:\alpha\in\mathbb{S}^1,\,j=1,2,3\}$, we observe that $D_2^*D_1F_3\nsim F_3$, and hence it easily follows that $\mathrm{Ad}_{F_3}(\Delta_3)\cap\mathrm{Ad}_{D_1^*D_2F_3}(\Delta_3)=\bbc$. This proves that $\Delta_3\cap\mathrm{Ad}_{u^*v}(\Delta_3)=\bbc$.\qed
\end{prf}

\begin{thm}\label{factoriality}
The von Neumann algebra $R_u \cap R_v$ is a {II}$_1$ subfactor of $R$ with $[R:R_u \cap R_v]=9$.
\end{thm}
\begin{prf}
By construction of the Hadamard subfactors $R_u\subset R\mbox{ and }R_v\subset R$, we have the following commuting squares
\begin{eqnarray}\label{zz}
\begin{matrix}
\mathrm{Ad}_{u_2}(\Delta_3\otimes M_3) &\subset & M_3\otimes M_3\\
\cup & &\cup\\
\mathrm{Ad}_u(\Delta_3) &\subset & \bbc\otimes M_3
\end{matrix}\quad\mbox{ and }\quad
\begin{matrix}
\mathrm{Ad}_{v_2}(\Delta_3\otimes M_3) &\subset & M_3\otimes M_3\\
\cup & &\cup\\
\mathrm{Ad}_v(\Delta_3) &\subset & \bbc\otimes M_3
\end{matrix}
\end{eqnarray}
Therefore, using \Cref{intersection} we obtain the commuting cube (see \Cref{commcube}) depicted in \Cref{commcubefac},
\begin{figure}
\centering
\begin{tikzpicture}
        \draw[thick,-<-=.5] (0,0)--(4,0);
        \draw[thick,->-=.7] (4,0)--(4,2);
        \draw[thick,->-=.5] (4,2)--(0,2);
        \draw[thick,-<-=.5] (0,2)--(0,0);
        \draw[thick,-<-=.5] (1,1)--(5,1);
        \draw[thick,->-=.5] (5,1)--(5,3);
        \draw[thick,->-=.5] (5,3)--(1,3);
        \draw[thick,-<-=.7] (1,3)--(1,1);
        \draw[thick,->-=.5] (0,0)--(1,1);
        \draw[thick,->-=.5] (4,0)--(5,1);
        \draw[->-=.5] (4,0)--(1,1);
        \draw[thick,->-=.5] (4,2)--(5,3);
        \draw[thick,->-=.5] (0,2)--(1,3);
        \draw[->-=.5] (4,2)--(1,3);

 \draw (-.3,-.2) node {\small{$\mathrm{Ad}_u(\Delta_3)$}};
        \draw (-1.2,2.1) node [scale =0.85]{{$\mathrm{Ad}_{u_2}(\Delta_3\otimes M_3)$}};
        \draw (5.7, 1.1) node [scale =0.85]{{$\mathrm{Ad}_v(\Delta_3)$}};
         \draw (.9, 1.2) node [scale =0.85]{{$\bbc\otimes M_3$}};
        \draw (.7,3.2) node [scale =0.85]{{$M_3\otimes M_3$}};
        \draw (6, 3.2) node [scale=0.85]{{$\mathrm{Ad}_{v_2}(\Delta_3\otimes M_3)$}};
        \draw (5.4,2.1) node [scale =0.85]{{$\mathrm{Ad}_{u_2W_2}(\bbc \otimes M_3)$}};
        \draw (4.25, -.1) node {\small{$\bbc$}};
\end{tikzpicture}
\caption{Commuting cube for proving factoriality}\label{commcubefac}
\end{figure}
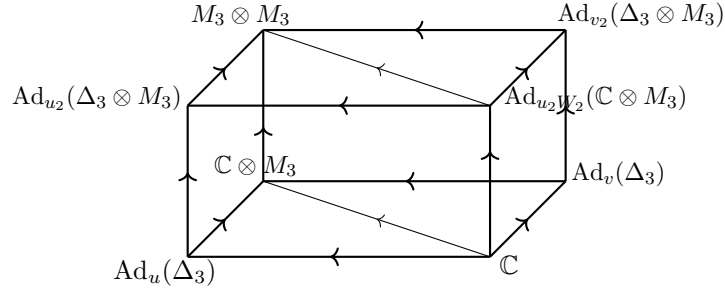
where the slice
\begin{eqnarray}\label{pcm}
\mathscr{I}:=\begin{matrix}
\mbox{Ad}_{u_2W_2}(\bbc \otimes M_3) &\subset & M_3\otimes M_3\\
\cup & &\cup\\
\bbc &\subset & \bbc\otimes M_3
\end{matrix}
\end{eqnarray}
is a commuting square due to \Cref{slice}. Observe that $\mathscr{I}$ is non-degenerate as norm of the inclusion matrices of both the horizontal embeddings are equal to $3$. The tower of basic construction for the inclusion $\bbc \otimes M_3 \subset M_3 \otimes M_3$ is given by, 
\begin{equation}\label{2.3}
\bbc \otimes M_3 \subset M_3 \otimes M_3 \subset^{f_1} M^{(2)}_3 \otimes M_3 \subset^{f_2} M^{(3)}_3 \otimes M_3 \subset \cdots  \subset^{f_{k-1}} M^{(k)}_3 \otimes M_3\subset\cdots
\end{equation}
where $f_k$'s are the Jones projection. We have the tower of basic construction for the construction of $R_u\subset R$ described by the following:
\begin{equation}\label{2.4}
\bbc \subset \Delta_3 \subset^{e_1} M_3 \subset^{e_2} \Delta_3 \otimes M_3 \subset^{e_1\otimes I_3} \Delta_3 \otimes M_3 \otimes M_3\subset^{e_2\otimes I_3}\cdots\subset  R
\end{equation}
For similar reason described in (Theorem $5.8$, \cite{BG1}) the limit of \Cref{2.3} is indeed $R$. That is, the hyperfinite type $II_1$ factor obtained in \Cref{2.3} is the same as that obtained in \Cref{2.4}. As $\mathscr{I}$ is a non-degenerate commuting square, the following tower of finite-dimensional $C^*$-algebras 
\begin{align*}
\bbc\subset L_0 \subset^{f_1} L_1\subset^{f_2} L_2 \subset\cdots\subset^{f_k} L_k\subset\cdots  
\end{align*}
is the Jones' tower of basic construction, where $L_0=\mbox{Ad}_{u_2W_2}(\bbc \otimes M_3)$ and $L_k=\{L_{k-1}, f_k\}^{\prime\prime}$ for $k \geq 1$.  Define  $R_{u,v}=\overline{\cup_{k}L_k}^{\,\mathrm{sot}}$. In other words, we obtain the basic construction of the non-degenerate commuting cube (see Section $3$, \cite{BG1}) depicted in \Cref{commcubefac}. By Corollary $5.7.4$ in \cite{JS}, it follows that $R_{u,v}$ is a $II_1$ factor and $[R:R_{u,v}]=9$ (since the norm of the inclusion matrix is $3$). By construction, $L_0\subset R_u\cap R_v$ and $f_k \in R_u \cap R_v$ for each $k$, and hence it follows that $R_{u,v}\subseteq R_u \cap R_v$. Now, $[R:R_{u,v}]=9$ implies that $[R_u:R_{u,v}]=3$ by the multiplicativity of the Jones index, as $[R:R_u]=3$. Therefore, $R_{u,v}\subset R_u$ is irreducible, which immediately gives factoriality of $R_u\cap R_v$, as $R_{u,v}\subseteq R_u\cap R_v$. Moreover, $[R_u\cap R_v:R_{u,v}]=1$, and hence $R_{u,v}=R_u\cap R_v$. Thus, we conclude that $R_u \cap R_v$ is a subfactor of $R$ with index $9$.\qed
\end{prf}

An immediate consequence is the following result. Since the proof is similar to Corollary $6.13$ in \cite{BG1}, we omit the details.
\begin{crlre}\label{refer later}
The quadruple $(R_u\cap R_v\subset R_u\,,R_v\subset R)$ of $II_1$ factors is obtained as an iterated basic construction of the non-degenerate commuting cube depicted in \Cref{commcubefac}.
\end{crlre}


\subsection{The relative commutant \texorpdfstring{$(R_u\cap R_v)^\prime\cap R$}{} and conjugacy}

In this subsection, we compute the relative commutant $(R_u\cap R_v)^\prime\cap R$ and establish the conjugacy of the Hadamard subfactors. Suppose that the following is a non-degenerate (symmetric) commuting square of connected inclusions of finite-dimensional $C^*$-algebras
\[
\begin{matrix}
A_{10} &\subset &A_{11} \cr \cup &\ &\cup\cr A_{00} &\subset & A_{01}
\end{matrix}
\]
Iterating the basic construction, we obtain the following ladder of non-degenerate commuting squares
\[
\begin{matrix}
A_{10} &\subset & A_{11} & \subset & A_{12} & \subset &  \cdots \cr
\rotatebox{90}{$\subset $} &\ &\rotatebox{90}{$\subset$} &\ & \rotatebox{90}{$\subset$} \cr
A_{00} &\subset & A_{01} & \subset & A_{02} & \subset & \cdots \,.
\end{matrix}
\]
Setting $A_{1,\infty}$ (resp.  $A_{0,\infty}$) as the GNS-completion of $\bigcup_k A_{1k}$ (resp. $\bigcup_k A_{0k}$), we obtain the hyperfinite subfactor $A_{0,\infty}\subset A_{1,\infty}$.

\begin{ppsn}[\cite{JS}]\label{ocneanucompactness}(Ocneanu compactness)
Let $A_{01},A_{10}, A_{0,\infty}$, and $A_{1,\infty}$ be as above. Then, 
$\big(A_{0,\infty}\big)^{\prime}\cap A_{1,\infty}=\big(A_{01})^{\prime}\cap A_{10}$.
\end{ppsn}

\begin{thm}\label{commutant}
The relative commutant $(R_u \cap R_v)^{'} \cap R$ is $\bbc\oplus\bbc\oplus\bbc$.
\end{thm}
\begin{prf}
Since $R_u \cap R_v=R_{u,v}$ by \Cref{factoriality}, using Ocneanu compactness we have the following
\begin{eqnarray}\label{pz}
(R_u \cap R_v)^{'} \cap R=R_{u,v}^{'}\cap R &=& (\mbox{Ad}_{u_2W_2}(\bbc \otimes M_3))^{'}\cap (\bbc \otimes M_3)\nonumber\\
&=& \mathrm{Ad}_{u_2W_2}(M_3 \otimes \bbc)\cap (\bbc \otimes M_3)
\end{eqnarray}
(recall $W_2$ from \Cref{notation}). For $x\in\mathrm{Ad}_{u_2 W_2}(M_3 \otimes \bbc)\cap (\bbc \otimes M_3)$, it is easy to verify that $\mathrm{Ad}_{W^*_2u^*_2}(x)=W^*_2(I_3\otimes F^*_3)D^*_{F_3}(I_3\otimes u^*xu) D_{F_3}(I_3\otimes F_3)W_2 \in M_3\otimes\bbc$ using \Cref{unitaries}. Now, define $y=u^*xu$ and observe the following
\[
W^*_2(F^*_3\otimes I_3)D^*_{{F}_{3}}(I_3\otimes y)D_{F_3}(F_3\otimes I_3)W_2=W^*_2\begin{bmatrix}
    Q_0(y) & Q_1(y) & Q_2(y)\\
    Q_2(y) & Q_0(y) & Q_1(y)\\
    Q_1(y) & Q_2(y) & Q_0(y)
\end{bmatrix}W_2
\]
(see \Cref{notation}). Since $\mathrm{Ad}_{W^*_2u^*_2}(x)\in M_3\otimes\bbc$, we get that $Q_{0}(y)=r_0I_3,\,Q_{1}(y)=r_1\sigma_2$ and $Q_2(y)=r_2\sigma_1$ for some $r_0, r_1, r_2 \in \bbc$. Therefore,
\[
y=u^*xu=\begin{bmatrix}
    r_0 & r_1 & r_2 \\
    r_2  & r_0 & r_1\\
    r_1 & r_2& r_0
\end{bmatrix}\,.
\]
From the matrix $u^*xu$ above, we conclude that $x\in\Delta_3$. So, we have $(R_u \cap R_v)^{'} \cap R\subseteq \Delta_3$ by \Cref{pz}. Conversely, note that $\mathrm{Ad}_{W^*_2u^*_2}(\Delta_3)\subseteq M_3 \otimes \bbc $. Therefore, $\Delta_3 \subseteq\mathrm{Ad}_{u_2W_2}( M_3  \otimes \bbc)\cap(\bbc\otimes M_3)=(R_u \cap R_v)^{'} \cap R$, which completes the proof.\qed
\end{prf}

\begin{lmma}\label{guru1}
For any $k\in\bbn$ and unitary matrix $w=DF_3\in M_3$ with $D\in\mathcal{U}(\Delta_3)$, we have $w_{2k}=(I_3^{(k)}\otimes D)(F_3)_{2k}$ (see \Cref{unitaries} for notations).
\end{lmma}
\begin{prf}
Follows by induction on $k$, together with the fact $D_w(w\otimes I_3)=D_{F_3}(F_3\otimes I_3)$.\qed
\end{prf}

\begin{thm}\label{R_u conjugate to R_v}
The pair of Hadamard subfactors arising from complex Hadamard matrices of order $3$ are conjugate to each other, and the conjugating unitary lies in the relative commutant of their intersection in the hyperfinite type $II_1$ factor $R$. 
\end{thm}
\begin{prf}
It is enough to show that for any $u=DF_3$, where $D$ is a diagonal unitary matrix in $M_3$, $R_u=\mbox{Ad}_D(R_{F_3})$. By \Cref{guru1}, we have $u_{2k}=D(F_3)_{2k}$ $\forall k\in\bbn$. This says that $B_{2k}^u=\mathrm{Ad}_D\big(B_{2k}^{F_3}\big)$ (see \Cref{not} in this regard) for all $k\in\bbn$. A simple limit argument shows that $R_u=\mathrm{Ad}_D(R_{F_3})$. Since $D\in\Delta_3\subseteq M_3$, by \Cref{commutant} the statement follows.\qed
\end{prf}

\begin{rmrk}\rm
\Cref{R_u conjugate to R_v} holds for the case of complex Hadamard matrices of order $2$ as well (\cite{BG1}); however, it fails for the case of Hadamard inequivalent complex Hadamard matrices of order $4$ (\cite{BG2}).
\end{rmrk}


\subsection{Characterization of \texorpdfstring{$R_u\cap R_v\subset R$}{} and vertex model}

In this subsection, we characterize the subfactor $R_u\cap R_v\subset R$ and draw its principal graph.
\begin{thm}\label{vertexintersection}
The subfactor $R_u\cap R_v \subset R$ is a vertex model subfactor of index $9$.
\end{thm}
\begin{prf}
By \Cref{factoriality}, the subfactor $R_u\cap R_v \subset R$ is obtained as iterated basic construction of the non-degenerate commuting square $\mathscr{I}$:
\[
\begin{matrix}
\mbox{Ad}_{u_2W_2}(\bbc \otimes M_3) &\subset & M_3\otimes M_3\\
\cup & &\cup\\
\bbc &\subset & \bbc\otimes M_3
\end{matrix}
\]
Observe that $\mbox{Ad}_{u_2W_2}(\bbc \otimes M_3)=\mbox{Ad}_{u_2W_2V_2}(M_3 \otimes\bbc)$, where $V_2:\sum_{i,j=1}^{3}E_{ij} \otimes E_{ji}$ is the flip operator. Therefore, $u_2W_2V_2$ is a bi-unitary matrix in $M_9$ (can be checked directly also), and consequently $R_u\cap R_v\subset R$ is a vertex model subfactor of index $9$.\qed
\end{prf}

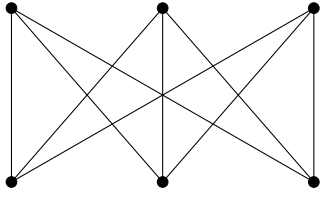
\begin{figure}
  \centering
  \begin{tikzpicture}
       \draw[] (0,0)--(0,2.3);
       \draw[] (0,2.3)--(-2,0);
       \draw[] (0,2.3)--(2,0);
       \draw[] (2,0)--(2,2.3);
       \draw[] (-2,0)--(-2,2.3);
       \draw[] (-2,2.3)--(0,0);
       \draw[] (2,2.3)--(0,0);
       \draw[] (-2,2.3)--(2,0);
       \draw[] (2,2.3)--(-2,0);
       \filldraw[black] (0,0) circle (2pt);
       \filldraw[black] (0,2.3) circle (2pt);
       \filldraw[black] (2,0) circle (2pt); 
       \filldraw[black] (2,2.3) circle (2pt);
      \filldraw[black] (-2,2.3) circle(2pt);
       \filldraw[black] (-2,0) circle (2pt);
 \end{tikzpicture}
    \caption{Principal graph of $R_u\cap R_v\subset R$}\label{PG}
 \end{figure}

\begin{thm}\label{ks type}
The subfactor $R_u\cap R_v \subset R$ is of depth $2$ and its principal graph is depicted in \Cref{PG}.
\end{thm}
\begin{prf}
We first claim that the bi-unitary matrix $u_2W_2V_2$ generating the vertex model subfactor $R_u\cap R_v \subset R$ in \Cref{vertexintersection} is equivalent (see \Cref{Sec 2} for the equivalence relation) to a bi-unitary permutation matrix. Indeed,
\begin{eqnarray*}
u_2W_2V_2 &=& (I_3\otimes u)D_{u}(u\otimes I_3)W_2V_2\nonumber\\
&=& (I_3\otimes u)D_{F_3}(F_3\otimes I_3)W_2V_2\nonumber\\
&=& (I_3\otimes u)(F^*_3 \otimes I_3)W_2\nonumber\\
&=& (F^*_3\otimes u)W_2.
\end{eqnarray*}
Since $W_2$ is a bi-unitary permutation matrix (see \Cref{notation}), our claim is justified. By Section $4$ in \cite{KSV}, we know that the vertex model subfactors corresponding to equivalent bi-unitary matrices are conjugate to each other. Since the bi-unitary matrix $u_2W_2V_2$ is equivalent to the bi-uintary permutation matrix $W_2$, by \Cref{bi-unitary} it follows that $\lambda=(\lambda_1, \lambda_2, \lambda_3)=(id, (123), (132))$ and $\rho=(\rho_1, \rho_2, \rho_3)=(id, id, id)$ are the corresponding elements in $P_n$ described in \Cref{Definition}. By \cite{KS}, we obtain the principal graph of the subfactor $R_u\cap R_v \subset R$ depicted in \Cref{PG}\footnote{The top left vertex in \Cref{PG} denotes $(R_u\cap R_v)^\prime\cap(R_u\cap R_v)=\bbc$ and the bottom three vertices denote the components corresponding to $(R_u\cap R_v)^\prime\cap R=\bbc^3$.}. It is now obvious that the subfactor $R_u\cap R_v \subset R$ is of depth $2$.\qed
\end{prf}

\newsection{Computation of a few key invariants}\label{Sec 5}

We explicitly compute a few key invariants for pairs of Hadamard subfactors $R_u,R_v\subset R$ each of index $3$ to understand their relative position. These include the Pimsner-Popa probabilistic number, the interior (and exterior) angle, the Sano-Watatani angle, and the Connes-St{\o}rmer relative entropy. As an application, we completely characterize when the quadruple $(R_u\cap R_v\subset R_u\,,\,R_v\subset R)$ of $II_1$ factors forms a commuting (and co-commuting) square.

\subsection{Pimsner-Popa probabilistic number}

For von Neumann subalgebras $\mathcal{P,Q}$ of a finite von Neumann algebra $\mathcal{M}$, the Pimsner-Popa probabilistic number \cite{PP} is defined by following
\begin{equation*}
\lambda(\mathcal{P,Q})=\text{sup}\{\lambda\geq 0\,:\,E_\mathcal{Q}^\mathcal{M}(x)\geq \lambda x\,\,\forall\,x\in {\mathcal{P}}_{+}\}.
\end{equation*}
In the case of a single subfactor $N\subset M,\,\lambda(M,N)=[M:N]^{-1}$ (with the convention $\frac{1}{0}=\infty$). For von Neumann algebras $\mathcal{P\subset M},\,\lambda(\mathcal{M,P})^{-1}$ is called the Pimsner-Popa index. For an irreducible quadruple of type $II_1$ factors $(N\subset P,Q\subset M)$ with $[M:N]<\infty,$ an explicit formula for $\lambda(P,Q)$ has been provided in \cite{B}. For a comprehensive detail, see Section $3$ in \cite{B} (see also Section $2.3$ in \cite{BG1}). However, no such formula for the non-irreducible situation is known. Our goal is to compute $\lambda(R_u,R_v)$ for the quadruple $(R_u\cap R_v\subset R_u, R_v\subset R)$ which is not irreducible (see \Cref{commutant}).

Recall from \Cref{Sec 4} that $u=D_1F_3$ and  $v=D_2F_3$, where $D_1=\mbox{diag}\{1,e^{i\alpha_1},e^{i\alpha_2}\}$ and $D_2=\mbox{diag}\{1,e^{i\beta_1},e^{i\beta_2}\}$ with $\alpha_1,\alpha_2,\beta_1,\beta_2\in[0,2\pi)$. First we need the following facts from \cite{BG1}.

\begin{ppsn}[Proposition $2.4$, \cite{BG1}]\label{popaadaptation}
\begin{enumerate}[$(i)$]
\item Let $\{M_n\}, \{A_n\}$ and $\{B_n\}$ be increasing sequences of von Neumann subalgebras of a finite von Neumann algebra $M$ such that $\{A_n\}, \{B_n\}\subset M$ and $M=\big(\bigcup_{n=1}^{\infty}M_n\big)^{\dprime}.$ If $A=\big(\bigcup_{n=1}^{\infty} A_n\big)^{\dprime}$ and $B=\big(\bigcup_{n=1}^{\infty} B_n\big)^{\dprime}$, then $\lambda(B,A)\geq \displaystyle \limsup\,\lambda(B_n,A_n).$
\item If in addition, $E_{A_{n+1}}E_{M_n}=E_{A_n}$ and $E_{B_{n+1}}E_{M_n}=E_{B_n}$ for $n\in \mathbb{N}$, then $\lambda(B,A)=\displaystyle \lim\,\lambda(B_n,A_n)$ decreasingly.
\end{enumerate}
\end{ppsn}

\begin{dfn}\label{hamdef}
Given a nonzero vector $\mathbf{w}\in\bbc^n$, the Hamming number is given by,
\[
h(\mathbf{w}):=\mbox{number of non-zero entries in }\mathbf{w}\,.
\]
\end{dfn}

\begin{thm}[Theorem $4.7$, \cite{BG1}]\label{finitemains}
If $\Delta_n$ and $U\Delta_n U^*$ are two Masas in $M_n$, where $U$ is a unitary matrix, then the Pimsner-Popa probabilistic number between them is given by the following
\[
\lambda(\Delta_n\,,\,U\Delta_n U^*)=\min_{1\leq i\leq n}\,\big(\mbox{$h\left(U^*\right)_i$}\big)^{-1}
\]
where $\left(U^*\right)_i$ is the $i$-th column of $U^*$.
\end{thm}

The above theorem is the backbone in proving \Cref{subrelation} (see Section $4$, \cite{BG1} for detail).

\begin{lmma}\label{rr}
We have $\lambda(u\Delta_3u^*,v\Delta_3v^*)=1/3$.
\end{lmma}
\begin{prf}
Note that $\lambda(u\Delta_3u^*,v\Delta_3v^*)=\lambda(\Delta_3,u^*v\Delta_3v^*u)$. Using the fact that for two complex numbers $z_1,z_2\in\mathbb{S}^1$, if $1+z_1+z_2=0$, then $(z_1,z_2)\in\{(\omega,\omega^2),(\omega^2,\omega)\}$ where $\omega$ is a primitive cube root of unity, we can conclude that all the entries of $u^*v$ are non-zero, where $u=D_1F_3\mbox{ and }v=D_2F_3$. Hence, by \Cref{finitemains} it follows that $\lambda(u\Delta_3u^*,v\Delta_3v^*)=1/3$.\qed
\end{prf}

\begin{lmma}\label{ss}
For von Neumann subalgebras $\mathcal{P,Q}$ of a finite von Neumann algebra $\mathcal{M}$, we have $\lambda(\mathcal{P,Q})\geq\lambda(\mathcal{M,Q})$.
\end{lmma}
\begin{prf}
Follows directly from their respective definitions.\qed
\end{prf}

\begin{thm}\label{PP final}
For the Hadamard subfactors $R_u\subset R\mbox{ and }R_v\subset R$, the Pimsner-Popa probabilistic number $\lambda(R_u,R_v)$ is equal to $1/3$.
\end{thm}
\begin{prf}
Recall the tower of basic construction for $R_u\subset R$ (similarly, $v$ in place of $u$) depicted in \Cref{main fig}. For any $k\in\bbn$, using \Cref{not} and \Cref{ss}, we get the following
\begin{eqnarray}\label{p0}
\lambda\big(B^u_{2k+1}\,,\,B^v_{2k+1}\big) &\geq& \lambda\big(A_{2k+1}\,,\,B^v_{2k+1}\big)\nonumber\\
&=& \lambda\big(\Delta_3\otimes M_3 \otimes M^{(k)}_3\,,\,v_{2k+1}(M_3 \otimes M^{(k)}_3) v^*_{2k+1}\big)\nonumber\\
&=& \lambda\big(\Delta_3 \otimes M_3 \otimes M^{(k)}_3\,,\,M_3 \otimes M^{(k)}_3\big)\,.
\end{eqnarray}
Here the last equality follows from the fact that $v_{2k+1}\in\Delta_3 \otimes M_3 \otimes M^{(k)}_3$. Now, notice that the trace on $A_{2k+1}=\Delta_3\otimes M_3 \otimes M^{(k)}_3$ is implemented by the restriction of the unique normalized trace on the type $I$ factor $A_{2k+2}=M_3^{(k+2)}$. A direct application of Theorem $6.1$ in \cite{PP} gives us $\lambda\big(\Delta_3\otimes M_3 \otimes M^{(k)}_3\,,\,M_3\otimes M^{(k)}_3\big)=1/3$ for any $k\geq 0$. Therefore, we have
\begin{eqnarray}\label{p1}
\lambda\big(B^u_{2k+1}\,,\,B^v_{2k+1}\big) &\geq& 1/3\,.
\end{eqnarray}
Since $\lambda(R_u,R_v)$ is limit of a decreasing sequence, by \Cref{popaadaptation} we have the following
\begin{eqnarray}\label{p2}
\lambda\big(B^u_{2k+1}\,,\,B^v_{2k+1}\big) &\leq& \lambda\big(B^u_1\,,\,B^v_1\big)\leq\lambda(B^u_0,B^v_0)=\lambda(u\Delta_3u^*,v\Delta_3v^*)=1/3
\end{eqnarray}
due to \Cref{rr}. Combining \Cref{p1,p2}, we get the following
\[
\lambda\big(B^u_{2k+1}\,,\,B^v_{2k+1}\big)=1/3
\]
for all $k\in\bbn$. By \Cref{popaadaptation}, we have $\lambda(R_u,R_v)=\lim_{k\to\infty}\,\lambda\big(B^u_{2k+1}\,,\,B^v_{2k+1}\big)=1/3$.\qed
\end{prf}


\subsection{Interior and exterior angle}

Consider intermediate subfactors $P,Q$ of a finite-index subfactor $N\subset M$. The notion of interior angle $\alpha^N_M(P,Q)$ and exterior angle $\beta^N_M(P,Q)$ between $P$ and $Q$ has been introduced in \cite{BDLR}. This angle is crucially used to improve the existing upper bound for the cardinality of the lattice of intermediate subfactors, and thereby answering a question of Longo. For more on angle, the readers are invited to the recent works in \cite{BG, GS, BGJ}. In this section, we compute these angles for the quadruple $(R_u\cap R_v\subset R_u\,,R_v\subset R)$.

\begin{dfn}[\cite{BDLR}]\label{defnangle}
Let $P$ and $Q$ be two intermediate subfactors of a finite-index subfactor $N\subset M$. The interior angle $\alpha^N_M(P,Q)$ between $P$ and $Q$ is defined by
\[
\alpha^N_M(P, Q)=\cos^{-1}{\langle v_P,v_Q\rangle}_2\,,
\]
where $v_P:=\frac{e_P-e_N}{{\lVert e_P-e_N\rVert}_2}$ (and similarly $v_Q$),$\,{\langle x,y\rangle}_2:=tr(y^*x)$ and ${\lVert x\rVert}_2:=(tr(x^*x))^{1/2}$. The exterior angle between $P$ and $Q$ is defined by $\beta^N_M(P, Q)=\alpha^M_{M_1}(P_1, Q_1)$, where $P_1$ (resp., $Q_1$) denotes the basic construction of $P\subset M$ (resp., $Q\subset M$).
\end{dfn}

It is known that a quadruple $(N\subset P,Q\subset M)$ of $II_1$ factors such that $[M:N]<\infty$ is a commuting square (resp., co-commuting) if and only if $\alpha^N_M(P, Q)$ (resp., $\beta^N_M(P, Q)$) equals $90^\circ$.

\begin{thm}\label{chotoangle}
The interior and exterior angle for the quadruple $(R_u\cap R_v\subset R_u\,,R_v\subset R)$ of $II_1$ factors are equal and given by the following
\[
\cos\big(\alpha^{R_u \cap R_v}_R(R_u, R_v)\big)=\frac{1}{9}\,\big{|}e^{-i(\alpha_1-\beta_1)}+e^{i(\alpha_2-\beta_2)}+e^{i(\alpha_1-\beta_1)}e^{-i(\alpha_2-\beta_2)}\big{|}^{2}\,,
\]
where $u=\mathrm{diag}\{1,e^{i\alpha_1},e^{i\alpha_2}\}F_3$ and $v=\mathrm{diag}\{1,e^{i\beta_1},e^{i\beta_2}\}F_3$.
\end{thm}
\begin{prf}
Using \Cref{refer later}, we have the commuting cube $\mathscr{C}_\infty$ depicted in \Cref{commangle}.
\begin{figure}
  \centering
  \begin{tikzpicture}
        \draw[thick,-<-=.5] (0,0)--(4,0);
        \draw[thick,->-=.7] (4,0)--(4,2);
        \draw[thick,->-=.5] (4,2)--(0,2);
        \draw[thick,-<-=.5] (0,2)--(0,0);
        \draw[thick,-<-=.5] (1,1)--(5,1);
        \draw[thick,->-=.5] (5,1)--(5,3);
        \draw[thick,->-=.5] (5,3)--(1,3);
        \draw[thick,-<-=.7] (1,3)--(1,1);
        \draw[thick,->-=.5] (0,0)--(1,1);
        \draw[thick,->-=.5] (4,0)--(5,1);
        \draw[->-=.5] (4,0)--(1,1);
        \draw[thick,->-=.5] (4,2)--(5,3);
        \draw[thick,->-=.5] (0,2)--(1,3);
        \draw[->-=.5] (4,2)--(1,3);

 \draw (-.3,-.25) node {$\text{Ad}_u(\Delta_3)$};
        \draw (-.48,2.1) node {$R_u$};
        \draw (5.8, 1) node {$\text{Ad}_v(\Delta_3)$};
         \draw (1.4, 1.24) node{$M_3$};
        \draw (.7,3.2) node{$R$};
        \draw (5.25, 3.2) node {$R_v$};
        \draw (4.8,2) node {$R_u\cap R_v$};
        \draw (4.25, -.25) node {$\bbc$};
 \end{tikzpicture}
    \caption{Commuting cube $\mathscr{C}_\infty$}\label{commangle}
 \end{figure}
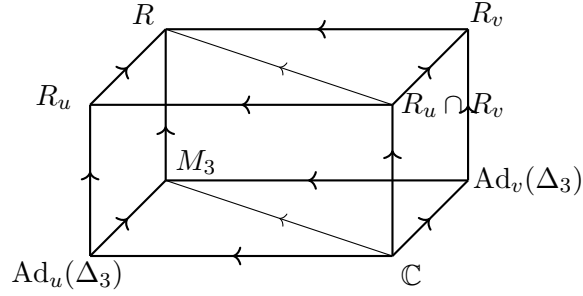
In particular, we have the commuting squares $(\bbc\subset \mathrm{Ad}_u(\Delta_3), R_u\cap R_v\subset R_u)$ and $(\bbc\subset \mathrm{Ad}_v(\Delta_3), R_u\cap R_v\subset R)$. Observe that $\{\lambda_1:= \text{Ad}_u(\sqrt{3}E_{11}), \lambda_2:= \text{Ad}_u(\sqrt{3}E_{22}), \lambda_3:= \text{Ad}_u(\sqrt{3}E_{33})$\} is a basis of $R_u/R_u\cap R_v$, and similarly  \{$ \mu_1:= \text{Ad}_v(\sqrt{3}E_{11}), \mu_2:= \text{Ad}_v(\sqrt{3}E_{22}), \mu_3:= \text{Ad}_v(\sqrt{3}E_{33})$\} is for $R_v/R_u\cap R_v$. Since, $(\bbc\subset M_3, R_u\cap R_v \subset R$) is a commuting square, and $\lambda_i, \mu_i \in M_3$ for $i=1,2,3$, we obtain $E^{R}_{R_u\cap R_v}(\lambda^*_i \mu_j)=E^{M_3}_{\bbc}(\lambda^*_i \mu_j)$ for all $ 1 \leq i,j \leq 3$. By the formula of the interior angle in \Cref{defnangle}, we observe the following
\begin{eqnarray}\label{zzx1}
\cos\big(\alpha^{R_u \cap R_v}_R(R_u, R_v)\big)&=&\frac{1}{2}\displaystyle{\sum_{i,j=1}^{3} tr(E^{R}_{R_u\cap R_v}(\lambda^*_i \mu_j)\mu^*_j \lambda_i)-\frac{1}{2}}\nonumber\\
&=& \frac{1}{2}\displaystyle {\sum_{i,j=1}^{3} tr(E^{M_3}_{\bbc}(\lambda^*_i\mu_j)\mu^*_j \lambda_i)-\frac{1}{2}}\nonumber\\
&=&\frac{1}{2}\displaystyle {\sum_{i,j=1}^{3} tr(\lambda^*_i \mu_j)tr(\mu^*_j \lambda_i)-\frac{1}{2}}.
\end{eqnarray}
A straightforward computation proceeding with the right-hand side in \Cref{zzx1} gives us the following~:
\begin{eqnarray}\label{zzx}
\cos\big(\alpha^{R_u \cap R_v}_R(R_u, R_v)\big) &=& \frac{1}{9}\,\big{|}e^{-i(\alpha_1-\beta_1)}+e^{i(\alpha_2-\beta_2)}+e^{i(\alpha_1-\beta_1)}e^{-i(\alpha_2-\beta_2)}\big{|}^{2}.
\end{eqnarray}
Since $(R_u \cap R_v \subset R_u, R_v\subset R)$ is a quadruple of $II_1$ factors with $[R:R_u\cap R_v]<\infty$, and $R_{u,v}=R_u\cap R_v \subset R$ is extremal with $[R_u: R_u \cap R_v]=[R : R_v]=3$, we have $\alpha^{R_u \cap R_v}_R=\beta^{R_u \cap R_v}_R$ by \cite{BDLR}, which concludes the proof.\qed
\end{prf}

An immediate consequence is the following.

\begin{crlre}\label{choto}
The quadruple $(R_u\cap R_v\subset R_u\,,R_v\subset R)$ of $II_1$ factors is a commuting square if and only if it is co-commuting square.
\end{crlre}

\begin{thm}[Characterization of a commuting square]\label{charac}
Let $u=\mathrm{diag}\{1,e^{i\alpha_1},e^{i\alpha_2}\}F_3$ and $v=\mathrm{diag}\{1,e^{i\beta_1},e^{i\beta_2}\}F_3$, with $\alpha_j,\beta_j\in[0,2\pi)$. The following are equivalent.
\begin{enumerate}[$(i)$]
\item The quadruple $(R_u\cap R_v \subset R_u, R_v \subset R)$ of $II_1$ factors is a commuting square;
\item The quadruple $(\bbc \subset u\Delta_3u^*, v\Delta_3v^* \subset M_3)$ of finite-dimensional $C^*$-algebras is a commuting square;
\item \[
\big(e^{i(\alpha_1-\beta_1)},e^{i(\alpha_2-\beta_2)}\big)\in\big\{(\omega,\omega),(\omega^2,\omega^2),(1,\omega),(1,\omega^2),(\omega,1),(\omega^2,1)\big\},
\]
where $\omega$ is a primitive cube root of unity.
\end{enumerate}
\end{thm}
\begin{prf}
\textbf{(i)$\implies$(ii):~} Consider the commuting cube depicted in \Cref{commangle}. By \Cref{klz}, if the roof $(R_u \cap R_v \subset R_u, R_v\subset R)$ is a commuting square, then the floor $(\bbc \subset u\Delta_3u^*, v\Delta_3v^* \subset M_3)$ is also a commuting square.
\smallskip

\textbf{(ii)$\implies$(iii):~} If $(\bbc \subset u\Delta_3u^*, v\Delta_3v^* \subset M_3)$ is a commuting square, then it is non-degenerate also as norm of the inclusion matrices of both the horizontal and vertical embeddings are equal. Hence, by \cite{JS} we get that $u^*v$ is a complex Hadamard matrix, where $u=\mathrm{diag}\{1,e^{i\alpha_1},e^{i\alpha_2}\}F_3$ and $v=\mathrm{diag}\{1,e^{i\beta_1},e^{i\beta_2}\}F_3$. Using the fact that the only solution of the equation $1+\xi+\eta=0$, where $|\xi|=|\eta|=1$, are given by $(\xi,\eta)\in\{(\omega,\omega^2),(\omega^2,\omega)\}$, where $\omega$ is a primitive cube root of unity, it is now a straightforward verification that the condition $(iii)$ is satisfied.
\smallskip

\textbf{(iii)$\implies$(i):~} Consider the following equation
\begin{eqnarray}\label{cm}
\overline{z_1}+z_2+z_1\overline{z_2}=0\,,
\end{eqnarray}
where $z_1,z_2\in\bbc$. Let $\omega$ be a primitive cube root of unity. If $(z_1,z_2)$ lies in the set in condition $(iii)$, then it is clear that $(z_1,z_2)$ satisfies \Cref{cm}. Put $z_1=e^{i(\alpha_1-\beta_1)}$ and $z_2=e^{i(\alpha_2-\beta_2)}$. Thus, the condition $(iii)$ implies that $\text{cos}(\alpha^{R_u \cap R_v}_R)=0$ by \Cref{zzx}. Therefore, the quadruple $(R_u \cap R_v \subset R_u, R_v\subset R)$ is a commuting square (and consequently, co-commuting square by \Cref{choto}).
\smallskip

This completes the proof.\qed
\end{prf}

\begin{rmrk}\rm
The quadruple ($R_u \cap R_v \subset R_u, R_v\subset R$) is always non-degenerate (similar to the $2\times 2$ case in \cite{BG1}), that is, $\overline{R_uR_v}^{\,\text{sot}}=R=\overline{R_vR_u}^{\,\text{sot}}$. Therefore, for suitable choice of $u\mbox{ and }v$, the quadruple $(R_u \cap R_v \subset R_u, R_v\subset R)$ provides another concrete example, apart from the one in \cite{BG1}, of a non-degenerate quadruple of $II_1$ factors that is neither commuting nor co-commuting square (see Theorem $7.1$ and Corollary $7.1$ in \cite{SW} and Theorem $3.21$ in \cite{GJ} in this regard).
\end{rmrk}


\subsection{Sano-Watatani angle operator}

To measure how far a quadruple is from being a commuting square, Sano-Watatani \cite{SW} introduced the notion of the `angle operator'. Let $\mathcal{H}$ be a Hilbert space and $\mathcal{K},\, \mathcal{L}$ be two different (closed) subspaces. Recall the angle operator $\Theta(p,q)$, where $p$ (resp. $q$) is the orthogonal projection onto $\mathcal{K}$ (resp. $\mathcal{L}$). The set $\text{Ang}(p,q)$ of angles between $p$ and $q$ is the subset of $[0,\pi/2]$ defined by the following (see Definition $2.1$ in \cite{SW}),
\begin{equation} 
\mathrm{Ang}(p,q)=\begin{cases}
\text{sp}~\Theta(p,q), & \text{if}~~ pq\neq qp.\\
\{\pi/2\},& \text{otherwise.}
\end{cases}
\end{equation}

Note that $\Theta$ is a positive operator and the spectrum of $\Theta$ is contained in $[0,\frac{\pi}{2}]$, but $0\mbox{ and }\frac{\pi}{2}$ are not eigenvalues.
\begin{dfn}[\cite{SW}]
Let $\mathcal{M}$ be a finite von Neumann algebra with a faithful normal tracial state $\mathrm{tr}$ and $\mathcal{P,Q}$ be von Neumann subalgebras of $\mathcal{M}$. The trace $\mathrm{tr}$ determines the normal faithful conditional expectations $E^{\mathcal{M}}_{\mathcal{P}}:\mathcal{M}\to\mathcal{P}$ and $E^{\mathcal{M}}_{\mathcal{Q}}:\mathcal{M}\to\mathcal{Q}$. They extend to the orthogonal projections $e_{\mathcal{P}}\mbox{ and }e_{\mathcal{Q}}$ on the GNS Hilbert space $L^2(\mathcal{M})$. The angle $\mathrm{Ang}_{\mathcal{M}}(\mathcal{P,Q})$ between $\mathcal{P}$ and $\mathcal{Q}$  is defined as follows~:
\begin{center}${\mathrm{Ang}}_{\mathcal{M}}(\mathcal{P,Q}):=\mathrm{Ang}(e_{\mathcal{P}},e_{\mathcal{Q}}).$\end{center}
\end{dfn}

For applications of the Sano-Watatani angle, see \cite{JX, GJ} for instance. Recall that \Cref{charac} gives us a complete characterization of the quadruple $(R_u\cap R_v \subset R_u, R_v \subset R)$ being a commuting square. When this is a commuting square, the Sano-Watatani angle is the singleton set $\{\pi/2\}$ by definition \cite{SW}. We compute the angle operator for the remaining cases below. For $u=\mathrm{diag}\{1,e^{i\alpha_1},e^{i\alpha_2}\}F_3$ and $v=\mathrm{diag}\{1,e^{i\beta_1},e^{i\beta_2}\}F_3$, consider the following complex number
\begin{eqnarray}\label{zxc}
\zeta=\frac{1}{3}\big(e^{i(\alpha_1-\beta_1)}+e^{-i(\alpha_2-\beta_2)}+e^{-i(\alpha_1-\beta_1)}e^{i(\alpha_2-\beta_2)}\big)\,.
\end{eqnarray}

Recall from \Cref{notation} the matrices $\sigma_1=E_{21}+E_{32}+E_{13}\mbox{ and }\sigma_2=E_{31}+E_{12}+E_{23}$ that are the representations of the alternating subgroup $A_3=\{\mbox{id},(123),(132)\}\unlhd S_3$ of even permutations, and also the generalized permutation matrix $Q_2(x)=x_{21}E_{21}+x_{32}E_{32}+x_{13}E_{13}$ for $x=(x_{ij})\in M_3$. The following technical lemma is crucial in the Sano-Watatani angle computation.

\begin{lmma}\label{to refer next}
Let $u=D_1F_3$ and $v=D_2F_3$, where $D_1,D_2$ are diagonal unitary matrices and $F_3$ is the $3\times 3$ Fourier matrix. The following identity holds.
\begin{enumerate}[$(i)$]
\item $E^{M_3}_{u\Delta_3u^*}(\mathrm{Ad}_{D_1}(\sigma_1))=\mathrm{Ad}_{D_1}(\sigma_1)$;
\item $E^{M_3}_{u\Delta_3u^*}(Q_2(x))=k_1\mathrm{Ad}_{D_1}(\sigma_1)$;
\item $E^{M_3}_{v\Delta_3v^*}(\mathrm{Ad}_{D_1}(\sigma_1))=k_2\mathrm{Ad}_{D_2}(\sigma_1)$;
\item $E^{M_3}_{u\Delta_3u^*}(\mathrm{Ad}_{D_2}(\sigma_1))=k_3\mathrm{Ad}_{D_1}(\sigma_1)$;
\end{enumerate}
where $k_1,k_2,k_3$ are certain constants.
\end{lmma}
\begin{prf}
These are straightforward verifications.\qed
\end{prf}

\begin{thm}\label{angle of sw}
Let $u=\mathrm{diag}\{1,e^{i\alpha_1},e^{i\alpha_2}\}F_3$ and $v=\mathrm{diag}\{1,e^{i\beta_1},e^{i\beta_2}\}F_3$. If the quadruple $(R_u\cap R_v \subset R_u, R_v \subset R)$ of $II_1$ factors is not a commuting square, the Sano-Watatani angle between the subfactors $R_u\mbox{ and }R_v$ is the singleton set $\{\arccos|\zeta|\}$, where $\zeta$ is the complex number defined in \Cref{zxc}.
\end{thm}
\begin{prf}
By \Cref{charac}, if the quadruple $(R_u\cap R_v \subset R_u, R_v \subset R)$ fails to become a commuting square, then so does the quadruple $(\bbc \subset u\Delta_3u^*, v\Delta_3v^* \subset M_3)$. Hence, the following positive matrix
\[
\widetilde{E}:=E^{A_0}_{B^u_0}E^{A_0}_{B^v_0}E^{A_0}_{B^u_0}-E^{A_0}_{C_0}=E^{M_3}_{\mathrm{Ad}_u(\Delta_3)}E^{M_3}_{\mathrm{Ad}_v(\Delta_3)}E^{M_3}_{\mathrm{Ad}_u(\Delta_3)}-E^{M_3}_\bbc
\]
(see \Cref{not} in this regard) is non-zero. Note that $E^{M_3}_\bbc$ is nothing but the unique normalized trace on $M_3$. By \Cref{notation}, we see that any $x\in M_3$ can be written as $x=Q_0(x)+Q_1(x)+Q_2(x)$.
\medskip

\noindent\textbf{Step 1:} To show that $\widetilde{E}(Q_0(x))=0$.

Since $u$ is a complex Hadamard matrix, so is its adjoint $u^*$. Therefore, we have the following commuting square
\begin{eqnarray}\label{cccc}
\begin{matrix}
\Delta_3 &\subset & M_3\\
\cup & & \cup\\
\bbc &\subset & \mathrm{Ad}_{u^*}(\Delta_3)
\end{matrix}
\end{eqnarray}
Now, $E^{M_3}_{\mathrm{Ad}_u(\Delta_3)}(Q_0(x))=\mathrm{Ad}_u\circ E^{M_3}_{\Delta_3}\circ\mathrm{Ad}_{u^*}(Q_0(x))$. Since $\mathrm{Ad}_{u^*}(Q_0(x))\in\mathrm{Ad}_{u^*}(\Delta_3)$ and the quadruple in \Cref{cccc} is a commuting square, we have $E^{M_3}_{\Delta_3}\circ\mathrm{Ad}_{u^*}(Q_0(x))=\mbox{tr}(Q_0(x))$. This immediately gives us that $E^{M_3}_{\mathrm{Ad}_u(\Delta_3)}E^{M_3}_{\mathrm{Ad}_v(\Delta_3)}E^{M_3}_{\mathrm{Ad}_u(\Delta_3)}(Q_0(x))=\mbox{tr}(Q_0(x))$, and consequently $\widetilde{E}(Q_0(x))=0$, which completes Step $1$.
\medskip

\noindent\textbf{Step 2:} To show that $\widetilde{E}^{\,2}(Q_2(x))=|\zeta|^2\widetilde{E}(Q_2(x))$ (where $\zeta$ is as in \Cref{zxc}).

Repeated application of \Cref{to refer next} shows the following
\begin{eqnarray*}
\widetilde{E}^{\,2}(Q_2(x)) &=& \widetilde{E}\big(E_{u\Delta_3u^*}E_{v\Delta_3v^*}E_{u\Delta_3u^*}(Q_2(x))\big)\nonumber\\
&=& \widetilde{E}\big(k_1k_2k_3\,\mathrm{Ad}_{D_1}(\sigma_1)\big)\nonumber\\
&=& k_1k_2^2k_3^2\,\mathrm{Ad}_{D_1}(\sigma_1)\nonumber\\
&=& k_2k_3\,\widetilde{E}(Q_2(x))\,.
\end{eqnarray*}
It remains to show that $k_2k_3=|\zeta|^2$. However, this is a direct verification by putting $D_1=\mathrm{diag}\{1,e^{i\alpha_1},e^{i\alpha_2}\}$ and $D_2=\mathrm{diag}\{1,e^{i\beta_1},e^{i\beta_2}\}$ in \Cref{to refer next}, and observing that $k_2=\zeta\mbox{ and }k_3=\overline{\zeta}$. This completes Step $2$.
\medskip

\noindent\textbf{Step 3:} To show that $\widetilde{E}^{\,2}(Q_1(x))=|\zeta|^2\widetilde{E}(Q_1(x))$ (where $\zeta$ is as in \Cref{zxc}).

Observe that for $x\in M_3$, there is $y\in M_3$ such that $Q_1(x)=Q_2(y)^*$. Then, we get the following chain of equalities
\[
\widetilde{E}^2(Q_1(x))=\widetilde{E}^2(Q_2(y)^*)=\widetilde{E}\big((\widetilde{E}(Q_2(y)))^*\big)=\big(\widetilde{E}^2(Q_2(y))\big)^*=|\zeta|^2\big(\widetilde{E}(Q_2(y))\big)^*
\]
by using Step $2$ at the last. Thus, we get the following desired equality
\[
\widetilde{E}^2(Q_1(x))=|\zeta|^2\widetilde{E}(Q_2(y))^*=|\zeta|^2\widetilde{E}(Q_2(y)^*)=|\zeta|^2\widetilde{E}(Q_1(x))\,,
\]
which completes Step $3$.
\medskip

Combining the above steps, in view of the fact that any $x\in M_3$ can be written as $x=Q_0(x)+Q_1(x)+Q_2(x)$, we get the following
\begin{eqnarray*}
\widetilde{E}^{\,2}(x) &=& \widetilde{E}\big(\widetilde{E}(Q_0(x))+\widetilde{E}(Q_1(x))+\widetilde{E}(Q_2(x))\big)\nonumber\\
&=& \widetilde{E}^{\,2}(Q_1(x)+Q_2(x))\nonumber\\
&=& |\zeta|^2\big(\widetilde{E}(Q_1(x))+\widetilde{E}(Q_2(x))\big)\\
&=& |\zeta|^2\widetilde{E}(x)\,.
\end{eqnarray*}
Thus, the following identity
\[
\big(E^{A_0}_{B^u_0}E^{A_0}_{B^v_0}E^{A_0}_{B^u_0}-E^{A_0}_{C_0}\big)^2=|\zeta|^2\big(E^{A_0}_{B^u_0}E^{A_0}_{B^v_0}E^{A_0}_{B^u_0}-E^{A_0}_{C_0}\big)
\]
holds. By \Cref{refer later} and \Cref{gggg}, the result now follows.\qed
\end{prf}


\subsection{Connes-St{\o}rmer relative entropy}

Generalizing the classical notion of conditional entropy in ergodic theory, Connes and Størmer \cite{CS} defined relative entropy $H(\mathcal{P|Q})$ between a pair of finite-dimensional von Neumann-subalgebras $\mathcal{P}\mbox{ and }\mathcal{Q}$ of a finite von Neumann algebra $\mathcal{M}$ equipped with a fixed faithful normal trace. Using the relative entropy as the main technical tool, a non-commutative version of the Kolmogorov-Sinai theorem is proved. Pimsner and Popa \cite{PP} discovered a surprising connection between relative entropy and the Jones index. They observed that the definition of the Connes-St{\o}rmer relative entropy does not depend on $\mathcal{P,Q}$ being finite-dimensional, so that one may also consider the relative entropy $H(\mathcal{P|Q})$ for arbitrary von Neumann subalgebras $\mathcal{P,Q}\subset\mathcal{M}$. For more on relative entropy, visit \cite{NS}.

\begin{dfn}[\cite{CS}]
Let $(\mathcal{M},\tau)$ be a finite von Neumann algebra and $\mathcal{P,Q\subseteq M}$ are von Neumann subalgebras. Let
\begin{align*}
\gamma &= \{x_j\in \mathcal{M}_+:\sum x_j=1,\,j=1,\ldots,n\}\,\,\,\mbox{ be a finite partition of unity},\cr
\eta &: [0,\infty)\longrightarrow\mathbb{R}\,\mbox{ be the continuous function }\,\,t\longmapsto -t\log t,\cr
H_\gamma(\mathcal{P|Q}) &:= \sum_{j=1}^n\left(\tau\circ\eta\,E_{\mathcal{Q}}^{\mathcal{M}}(x_j)-\tau\circ\eta\,E_{\mathcal{P}}^{\mathcal{M}}(x_j)\right).
\end{align*}
Then, $H(\mathcal{P|Q}):=\sup_{\gamma}\,H_\gamma(\mathcal{P|Q})$ is the Connes-St{\o}rmer relative entropy between $\mathcal{P}\mbox{ and }\mathcal{Q}$.
\end{dfn}

\begin{dfn}[\cite{choda},\cite{choda2}]
Suppose that $(\mathcal{M},\tau)$ be a finite von Neumann algebra and $\mathcal{P,Q\subseteq M}$ are von Neumann subalgebras. Let
\begin{align*}
\gamma &= \{x_j\in\mathcal{P}_+:\sum x_j=1,\,j=1,\ldots,n\}\,\,\,\mbox{ be a finite partition of unity},\cr
\eta &: [0,\infty)\longrightarrow\mathbb{R}\,\mbox{ be the continuous function }\,\,t\longmapsto -t\log t,\cr
h_\gamma(\mathcal{P|Q}) &:= \sum_{j=1}^n\left(\tau\circ\eta\,E_{\mathcal{Q}}^{\mathcal{M}}(x_j)-\tau\circ\eta\,(x_j)\right).
\end{align*}
Then, $h(\mathcal{P|Q}):=\sup_{\gamma}\,h_\gamma(\mathcal{P|Q})$ is called the modified Connes-St{\o}rmer relative entropy between $\mathcal{P}\mbox{ and }\mathcal{Q}$.
\end{dfn}

If $\mathcal{M}$ is abelian, then $H(\mathcal{P|Q})=h(\mathcal{P|Q})$. Thus, $h$ also generalizes the classical relative entropy. Moreover, if $\mathcal{P}\subset\mathcal{M}$ then $H(\mathcal{M|P})=h(\mathcal{M|P})$. Also, it is known that $0\leq h(\mathcal{P|Q})\leq H(\mathcal{P|Q})$. For the case of commuting square, these two relative entropies agree \cite{choda2}.

\begin{ppsn}[Proposition $2.6$, \cite{BG1}]\label{popaadaptation2}
\begin{enumerate}[$(i)$]
\item Let $\{M_n\}, \{A_n\}$ and $\{B_n\}$ be increasing sequences of von Neumann subalgebras of a finite von Neumann algebra $M$ such that $\{A_n\}, \{B_n\}\subset M$ and $M=\big(\bigcup_{n=1}^{\infty}M_n\big)^{\dprime}.$ If $A=\big(\bigcup_{n=1}^{\infty} A_n\big)^{\dprime}$ and $B=\big(\bigcup_{n=1}^{\infty} B_n\big)^{\dprime}$, then $H(B|A)\leq \displaystyle \liminf\,H(B_n|A_n).$
\item If in addition, $E_{A_{n+1}}E_{M_n}=E_{A_n}$ and $E_{B_{n+1}}E_{M_n}=E_{B_n}$ for $n\in \mathbb{N}$, then $H(B|A)=\displaystyle \lim\,H(B_n|A_n)$ increasingly.
\end{enumerate}
\end{ppsn}

Similar statement also holds for $h$ in place of $H$.

\begin{ppsn}\label{choda setup}
Let $R_{F_3}\subset R$ be the Hadamard subfactor corresponding to the Fourier matrix $F_3$. Then, $R=R_{F_3}\rtimes^\theta\bbz_3$ for the outer action $\theta$ defined by $\theta_g(x)=\mathrm{Ad}_{\mathrm{diag}\{1,\omega,\omega^2\}}(x)$ for $x\in R_{F_3}$, where $g$ is a generator of the finite cyclic group $\bbz_3$.
\end{ppsn}
\begin{prf}
First observe that $\mbox{diag}\{1,\omega,\omega^2\}F_3=F_3\sigma_1$, where $\sigma_1=E_{13}+E_{21}+E_{32}$. This says that $\mathrm{Ad}_{\mathrm{diag}\{1,\omega,\omega^2\}}(u\Delta_3u^*)\subseteq u\Delta_3u^*$ because $\sigma_1$ is a permutation matrix, in which case it normalizes $\Delta_3$. Since $\mathrm{diag}\{1,\omega,\omega^2\}$ commutes with all the Jones' projections in the tower of basic construction for $R_{F_3}\subset R$, a simple induction argument shows that $\mathrm{Ad}_{\mathrm{diag}\{1,\omega,\omega^2\}}(B^u_{2k})\subseteq B^u_{2k}$ and $\mathrm{Ad}_{\mathrm{diag}\{1,\omega,\omega^2\}}(B^u_{2k+1})\subseteq B^u_{2k+1}$ for all $k\in\bbn$ (recall \Cref{not} in this regard). This establishes that the action $\theta$ is well-defined. Now, by the construction of the Hadamard subfactor $R_{F_3}\subset R$, we have the following commuting square
\[
\begin{matrix}
R_{F_3} &\subset & R\\
\cup & & \cup\\
\bbc &\subset & \Delta_3
\end{matrix}
\]
Since $\mbox{diag}\{1,\omega,\omega^2\}\in\Delta_3$ and $\Delta_3\cap R_{F_3}=\bbc$, it follows that $\theta$ is outer with outer period $3$. We then have $R_{F_3}\subsetneq R_{F_3}\rtimes^\theta\bbz_3\subseteq R$. Since $R_{F_3}\subset R$ is irreducible (being a Hadamard subfactor), it follows that $R_{F_3}\rtimes^\theta\bbz_3$ is a $II_1$ factor. Moreover, $R_{F_3}\rtimes^\theta\bbz_3=R$ because $[R:R_{F_3}]=3$ and Jones index is multiplicative.\qed
\end{prf}

\begin{thm}\label{ent of cs}
Let $u=\mathrm{diag}\{1,e^{i\alpha_1},e^{i\alpha_2}\}F_3$ and $v=\mathrm{diag}\{1,e^{i\beta_1},e^{i\beta_2}\}F_3$ be $3\times 3$ complex Hadamard matrices. For the pair of Hadamard subfactors $R_u\subset R$ and $R_v\subset R$, we have the following:
\begin{enumerate}[$(i)$]
\item $H(R|R_u \cap R_v)=3 \log3$ and $H(R_u|R_u\cap R_v)=H(R_v | R_u \cap R_v)=\log 3$;
\item $h(R_u|R_v)=\eta\left(\frac{1}{9}\,|1+e^{i(\beta_1-\alpha_1)}+e^{i(\beta_2-\alpha_2)}|^2\right)+\eta\left(\frac{1}{9}\,|1+e^{i(\beta_1-\alpha_1)}\omega+e^{i(\beta_2-\alpha_2)}\omega^2|^2\right)\\
\hspace*{2cm} +\eta\left(\frac{1}{9}\,|1+e^{i(\beta_1-\alpha_1)}\omega^2+e^{i(\beta_2-\alpha_2)}\omega|^2\right),$\\
where $\omega$ is a primitive cube root of unity;
\item We have $h(R_u|R_v)\leq H(R_u|R_v) \leq  \log 3 $. Further, when the quadruple $(R_u\cap R_v \subset R_u, R_v \subset R)$ of $II_1$ factors is a commuting square, then $h(R_u|R_v)= H(R_u|R_v)=-\log\lambda(R_u,R_v)=\log 3$.
\end{enumerate}
\end{thm}
\begin{prf}
Since $R_u\cap R_v\subset R$ is a factor by \Cref{factoriality}, part $(i)$ follows from \cite{PP}. For the second part, first by Corollary $3$ in \cite{choda}, we have the following
\begin{eqnarray}\label{xxx}
h(u\Delta_3u^*|v\Delta_3 v^*)&=& h(\Delta_3|u^*v\Delta_3v^*u)\nonumber\\
&=& \frac{1}{3}\sum_{i,j=1}^{3}\eta(|(u^*v)_{ij}|^2)\nonumber\\
&=& \eta\left(\frac{1}{9}\,\big|1+e^{i(\beta_1-\alpha_1)}+e^{i(\beta_2-\alpha_2)}\big|^2\right)+\eta\left(\frac{1}{9}\,\big|1+e^{i(\beta_1-\alpha_1)}\omega+e^{i(\beta_2-\alpha_2)}\omega^2\big|^2\right)\nonumber\\
& & +\eta\left(\frac{1}{9}\,\big|1+e^{i(\beta_1-\alpha_1)}\omega^2+e^{i(\beta_2-\alpha_2)}\omega\big|^2\right),
\end{eqnarray}
where $u=D_1F_3$ and $v=D_2F_3$, with $D_1=\mathrm{diag}\{1,e^{i\alpha_1},e^{i\alpha_2}\}$ and $D_2=\mathrm{diag}\{1,e^{i\beta_1},e^{i\beta_2}\}$. By \Cref{R_u conjugate to R_v}, we have $R_u=\mbox{Ad}_{D_1}(R_{F_3})$ and $R_v=\mbox{Ad}_{D_2}(R_{F_3})$. Therefore, $h(R_u|R_v)=h\big(R_{F_3}|\mathrm{Ad}_{D_1^*D_2}(R_{F_3})\big)$. Since $D_1^*D_2=\mbox{diag}\{1,e^{i(\beta_1-\alpha_1)},e^{i(\beta_2-\alpha_2)}\}$, we see that $D_1^*D_2=\gamma_0I_3+\gamma_1\,\mbox{diag}\{1,\omega,\omega^2\}+\gamma_2\,\mbox{diag}\{1,\omega^2,\omega\}$, where
\begin{eqnarray}\label{ref}
& & \gamma_0=\frac{1}{3}\big(1+e^{i(\beta_1-\alpha_1)}+e^{i(\beta_2-\alpha_2)}\big),\nonumber\\
& & \gamma_1=\frac{1}{3}\big(1+e^{i(\beta_1-\alpha_1)}\omega^2+e^{i(\beta_2-\alpha_2)}\omega\big),\nonumber\\
& & \gamma_2=\frac{1}{3}\big(1+e^{i(\beta_1-\alpha_1)}\omega+e^{i(\beta_2-\alpha_2)}\omega^2\big)\,.
\end{eqnarray}
By \Cref{choda setup}, $R_{F_3}$ is a subfactor of $R_{F_3}\rtimes^\theta\bbz_3$, and to find the value of $h\big(R_{F_3}|\mathrm{Ad}_{D_1^*D_2}(R_{F_3})\big)$, we apply Theorem $3.14$ in \cite{choda2} (recall the action of $\theta$ here) to get the following inequality
\[
h(R_u|R_v)=h\big(R_{F_3}|\mathrm{Ad}_{D_1^*D_2}(R_{F_3})\big)\leq\sum_{j=0}^2\eta\tau(\gamma_j\gamma_j^*)=\sum_{j=0}^2\eta(|\gamma_j|^2),
\]
where $\gamma_j$'s are as defined in \Cref{ref}. Hence, by \Cref{xxx} we obtain the inequality $h(R_u|R_v)\leq h(u\Delta_3u^*|v\Delta_3 v^*)$. Since the reverse inequality $h(R_u|R_v)\geq h(u\Delta_3 u^*|v\Delta v^*)$ is obvious by \Cref{popaadaptation2}, part $(ii)$ is concluded.

Finally, $h(R_u|R_v)\leq H(R_u|R_v)$ follows from their definition, and $H(R_u|R_v)\leq H(R_u|R)+H(R|R_v)=\log[R:R_v]=\log 3$ follows from \cite{PP}. When the quadruple $(R_u\cap R_v \subset R_u, R_v \subset R)$ of $II_1$ factors is a commuting square, we get that $h(R_u|R_v)=H(R_u|R_v)$ by \cite{choda2}, and $H(R_u|R_v)=H(R_u|R_u\cap R_v)$ by \cite{WW}. Thus, in this case we have the following chain of equalities
\[
h(R_u|R_v)=H(R_u|R_v)=\log[R_u:R_u\cap R_v]=\log 3=-\log\lambda(R_u,R_v)
\]
using \Cref{factoriality} and \Cref{PP final} respectively.\qed
\end{prf}
\medskip

We conclude this article with the following question.
\smallskip

\noindent \textbf{Open Question:} What is the value of the Connes-St{\o}rmer relative entropy $H(R_u|R_v)$ in general?
\smallskip

A major difficulty in attacking this question is that even the value $H(u\Delta_3u^*|v\Delta_3v^*)$, which is the first term of the sequence of relative entropies (in view of \Cref{popaadaptation2}) in the tower of basic construction (depicted in \Cref{main fig}), seems unknown in the literature (see \cite{PAM}).
\bigskip

\section*{Acknowledgements}
Keshab Chandra Bakshi acknowledges the support of DST INSPIRE Faculty fellowship\\ DST/INSPIRE/04/2019/002754, and Satyajit Guin acknowledges the support of SERB grant\\ MTR/2021/000818.

\bigskip

\bigskip

\noindent {\em Department of Mathematics and Statistics},\\
{\em Indian Institute of Technology Kanpur},\\
{\em Uttar Pradesh $208016$, India}
\medskip

\noindent {Keshab Chandra Bakshi:} keshab@iitk.ac.in; bakshi209@gmail.com\\
{Satyajit Guin:} sguin@iitk.ac.in\\
{Guruprasad:} guruprasad8909@gmail.com


\bigskip

\begin{thebibliography}{10}

\bibitem{B}
Bakshi, K. C.: 
\newblock A short note on relative entropy for a pair of intermediate subfactors.
\newblock {\em Proc. Amer. Math. Soc.} 150 (2022), no. 9, 3899--3913.

\bibitem{BDLR}
Bakshi, K. C.; Das, S.; Liu, Z.; Ren, Y.:
\newblock An angle between intermediate subfactors and its rigidity.
\newblock {\em Trans. Amer. Math. Soc.} 371 (2019), no. 8, 5973--5991.

\bibitem{BG1}
Bakshi, K. C.; Guin, S.:
\newblock Relative position between a pair of spin model subfactors.
\newblock {\em J. Aust. Math. Soc.} \href{https://www.cambridge.org/core/journals/journal-of-the-australian-mathematical-society/article/relative-position-between-a-pair-of-spin-model-subfactors/0277143FCA0DEE91672D67B4493D2EA0}{DOI:10.1017/S1446788725000035}.

\bibitem{BG2}
Bakshi, K. C.; Guin, S.:
\newblock A family of subfactors arising from a pair of complex Hadamard matrices.
\newblock {\em Internat. J. Math.} 36 (2025), no. 2, Paper No. 2450071.

\bibitem{BGJ}
Bakshi, K. C.; Guin, S; Jana, D.:
\newblock A few remarks on intermediate subalgebras of an inclusion of $C^*$-algebras.
\newblock {\em Infin. Dimens. Anal. Quantum Probab. Relat. Top.} \href{https://www.worldscientific.com/doi/10.1142/S0219025725500067?srsltid=AfmBOooRmZ6f6UvrcMUkp6F9aBpIFRSQtbAnvv-_oXiKSmECA0RVMwXq}{DOI:10.1142/S0219025725500067}.

\bibitem{BG}
Bakshi, K. C.; Gupta, V. P.:
\newblock Lattice of intermediate subalgebras. 
\newblock {\em J. Lond. Math. Soc. $(2)$} 104 (2021), no. 5, 2082--2127.

\bibitem{BINA}
Bhattacharyya, B.: 
\newblock Group actions on graphs related to Krishnan-Sunder subfactors.
\newblock {\em Trans. Amer. Math. Soc.} 355 (2003), no. 2, 433--463.

\bibitem{choda}
Choda, M:
\newblock Relative entropy for maximal abelian subalgebras of matrices and the entropy of unistochastic matrices.
\newblock{\em Internat. J. Math.} 19 (2008), no. 7, 767--776.

\bibitem{choda2}
Choda, M:
\newblock Conjugate pairs of subfactors and entropy for automorphisms.
\newblock{\em Internat. J. Math.} 22 (2011), no. 4, 577--592.

\bibitem{CS}
Connes, A.; Størmer, E.:
\newblock Entropy for automorphisms of $II_1$ von Neumann algebras.
\newblock {\em Acta Math.} 134 (1975), no. 3-4, 289--306.

\bibitem{GHJ}
Goodman, F.; de la Harpe, P.; Jones, V. F. R.:
\newblock Coxeter graphs and towers of algebras.
\newblock Mathematical Sciences Research Institute Publications, 14. {\em Springer-Verlag, New York}, 1989.

\bibitem{GJ}
Grossman, P.; Jones, V. F. R.:
\newblock Intermediate subfactors with no extra structure.
\newblock {\em J. Amer. Math. Soc.} 20 (2007), no. 1, 219--265.

\bibitem{GS}
Gupta, V.P.; Sharma, D.: 
\newblock On possible values of the interior angle between intermediate subalgebras.
\newblock {\em J. Aust. Math. Soc.} 117 (2024), no. 1, 44--66.

\bibitem{Jo}
Jones, V. F. R.:
\newblock Index for subfactors.
\newblock {\em Invent. Math.} 72 (1983), no. 1, 1--25.

\bibitem{JS}
Jones, V. F. R.; Sunder, V. S.:
\newblock Introduction to subfactors.
\newblock London Mathematical Society Lecture Note Series, 234. {\em Cambridge University Press, Cambridge}, 1997.

\bibitem{JX}
Jones, V. F. R.; Xu, F.:
\newblock Intersections of finite families of finite index subfactors.
\newblock {\em Internat. J. Math.} 15 (2004), no. 7, 717--733.

\bibitem{Jo1}
Jones, V. F. R.:
\newblock Two subfactors and the algebraic decomposition of bimodules over $II_1$ factors.
\newblock {\em Acta Math. Vietnam.} 33 (2008), no. 3, 209--218.

\bibitem{Jo2}
Jones, V. F. R.:
\newblock Planar algebra, I.
\newblock{\em New Zealand J. Math.} 52, 1--107, 2021.

\bibitem{K}
Kawahigashi, Y.:
\newblock Classification of paragroup actions in subfactors.
\newblock {\em Publ. Res. Inst. Math. Sci.} 31 (1995), no. 3, 481--517.

\bibitem{KSV} 
Krishnan, U.; Sunder, V. S.;  Varughese, C.: 
\newblock On some subfactors of integer index arising from vertex models.
\newblock {\em J. Funct. Anal.} 140 (1996), no. 2, 449--471.

\bibitem{KS} 
Krishnan, U.; Sunder, V. S.: 
\newblock On biunitary permutation matrices and some subfactors of index $9$.
\newblock {\em Trans. Amer. Math. Soc.} 348 (1996), no. 12, 4691--4736.

\bibitem{NS}
Neshveyev, S; Størmer, E.:
\newblock Dynamical entropy in operator algebras.
\newblock Springer-Verlag, Berlin, 2006.

\bibitem{N}
Nicoara, R.:
\newblock Subfactors and Hadamard matrices.
\newblock {\em  J. Operator Theory} 64 (2010), no. 2, 453--468.

\bibitem{Oc}
Ocneanu, A.:
\newblock Operator Algebras, topology and subgroups of quantum symmetry--construction of subgroups of quantum groups.
\newblock {\em Taniguchi Conference on Mathematics Nara '98}, 235--263, Adv. Stud. Pure Math., 31, {\em Math. Soc. Japan, Tokyo}, 2001.

\bibitem{PAM}
Petz, D; András, S.; Mihály, W.:
\newblock Complementarity and the algebraic structure of four-level quantum systems.
\newblock {\em Infin. Dimens. Anal. Quantum Probab. Relat. Top.} 12 (2009), no. 1, 99--116.

\bibitem{PP}
Pimsner, M.; Popa, S.:
\newblock Entropy and index for subfactors.
\newblock {\em Ann. Sci. École Norm. Sup.} (4) 19 (1986), no. 1, 57--106.

\bibitem{SW}
Sano, T.; Watatani, Y.:
\newblock Angles between two subfactors.
\newblock {\em J. Operator Theory} 32 (1994), no. 2, 209--241.

\bibitem{WW}
Watatani, Y.; Wierzbicki, J.:
\newblock Commuting squares and relative entropy for two subfactors.
\newblock {\em J. Funct. Anal.} 133 (1995), no. 2, 329--341.

\end{thebibliography}
\end{document}